\documentclass [a4paper, 10pt] {article}

\usepackage{amssymb}
\usepackage{latexsym}
\usepackage[mathscr]{eucal}
\usepackage{theorem}
\usepackage[all]{xy}

\newenvironment{proof}{\par \noindent{\bf Proof: }}{\hspace{\stretch{1}} $\Box$ \par \mbox{}}
\newcommand{\noproof}{\hspace{\stretch{1}} $\Box$}

\newtheorem{theorem}{Theorem}[section]
\newtheorem{proposition}[theorem]{Proposition}
\newtheorem{lemma}[theorem]{Lemma}

{\theorembodyfont{\rmfamily}
\newtheorem{definition}[theorem]{Definition}
\newtheorem{example}[theorem]{Example}
}

\newenvironment{theorem*}{\par \medskip \noindent{\bf Theorem }}{\par \mbox{}}
\newenvironment{lemma*}{\par \medskip \noindent{\bf Theorem }}{\par \mbox{}}

\newcommand{\Hom}{\mathop{Hom}}

\newcommand{\Ob}{\mathop{Ob}}

\newcommand{\gtp}{\mathop{\hat{\otimes}}}

\newcommand{\Or}{\mathop{Or}}
\newcommand{\EG}{\underline{E}G}

\begin{document}

\title{$C^\ast$-categories, Groupoid Actions, Equivariant $KK$-theory,
  and the Baum-Connes Conjecture}

\author{Paul D. Mitchener}

\maketitle

\section*{Abstract}
In this article we give a characterisation of the Baum-Connes assembly
map with coefficients.  The technical tools needed are the $K$-theory
of $C^\ast$-categories, and equivariant $KK$-theory in the world of groupoids.

Keywords: $C^\star$-category, $KK$-theory, Assembly map

AMS 2000 Mathematics Subject Classification: 19K35, 55P42

\tableofcontents

\section{Introduction}

Let $G$ be a discrete group.  Any given $G$-homotopy-invariant functor from the
category of $G$-$CW$-complexes to the category of spectra has a
universal approximation by a generalised $G$-equivariant homology theory.  To be
specific, we have the following result, proved by Davis and
L{\"u}ck in \cite{DL}.

\begin{theorem} \label{1}
Let $\mathbb E$ be a $G$-homotopy-invariant functor from the category
of proper $G$-$CW$-complexes to the category of spectra.  Then there
is a $G$-homotopy-invariant excisive functor ${\mathbb E}^\%$ and a
natural transformation $\alpha \colon {\mathbb E}^\% \rightarrow
{\mathbb E}$ such that the map
$$\alpha \colon {\mathbb E}^\% (G/H) \rightarrow {\mathbb E}(G/H)$$
is a stable equivalence for every finite subgroup, $H$, of the group
$G$.

Further, the pair $({\mathbb E}^\% , \alpha )$ is unique up to weak
equivalence.
\noproof
\end{theorem}

Here, a functor ${\mathbb E}^\%$ is called {\em excisive} when the
collection of functors $X\mapsto \pi_\star {\mathbb E}^\% (X)$ form a $G$-equivariant generalised homology
theory.  The natural transformation $\alpha \colon {\mathbb E}^\% \rightarrow
{\mathbb E}$ is called the {\em assembly map} associated to the
functor $\mathbb E$.  

The above theorem, or rather a slight generalisation, can be
used to describe several standard maps that appear in isomorphism
conjectures.  For example, the map appearing in the Farell-Jones
isomorphism conjecture (see \cite{FJ}) readily fits into the framework
described by Davis and L{\"u}ck.

Now, let $A$ be a $G$-$C^\ast$-algebra.  For a proper $G$-space $X$, one can
define $G$-equivariant $K$-homology groups, $K_n^G (X;A)$, with
coefficients in the $G$-$C^\ast$-algebra $A$.  There is a canonical
map
$$\beta \colon K_n^G (X;A)\rightarrow K_n (A\rtimes_r G)$$

Here $A\rtimes_r G$ is the reduced crossed product of the
$C^\ast$-algebra $A$ with the group $G$.
The map $\beta$ is termed the {\em Baum-Connes assembly map}

A proper $G$-$CW$-complex $\EG$ is called a {\em classifying space for proper 
actions of $G$} if for a given subgroup $H\leq G$ the fixed point set 
$\EG^H$ is contractible when $H$ is finite, and empty when $H$ is
infinite.  The {\em Baum-Connes conjecture} with
{\em coefficients} in the $G$-$C^\ast$-algebra $A$ is the assertion
that the Baum-Connes assembly map
$$\beta \colon K_n^G (X;A)\rightarrow K_n (A\rtimes_r G)$$
is an isomorphism.

The reader is urged to consult \cite{BCH} for a full description of
the Baum-Connes conjecture and details of some of its geometric and
algebraic implications.

In this paper we give a description of the Baum-Connes assembly map at
the level of spectra that fits into the framework described by Davis
and L{\"u}ck.  In order to use theorem \ref{1}, we need to consider
groupoids.  Actions on spaces naturally lead to groupoids because of
the following standard construction.

\begin{definition}
Let $X$ be a $G$-space.  Then we write $\overline{X}$ to denote the
category in which the collection of objects is that set $X$, and the
morphism sets are defined by writing
$$\Hom (x,y)_{\overline{X}} = \{ g\in G \ |\ xg = y \}$$
\end{definition}

Every morphism in the category $\overline{X}$ is invertible, so the
category $\overline{X}$ is a groupoid.  In this paper we ignore the
topology of the space $X$ when considering the groupoid $\overline{X}$.

If $\cal G$ is a groupoid, we define a ${\cal G}$-$C^\ast$-algebra to
be a functor from the groupoid $\cal G$ to the category of
$C^\ast$-algebras.  If $A$ is a $\cal G$-$C^\ast$-algebra, there is a
natural notion of the {\em reduced crossed product}, $A\rtimes_r {\cal
G}$.  When $\cal G$ is a groupoid rather than a group, this reduced
crossed product is not a $C^\ast$-algebra, but rather a more general
object called a {\em $C^\ast$-category}, as defined in \cite{GLR}.  

If $f\colon {\cal G}\rightarrow {\cal H}$ is a faithful functor
between groupoids, and $A$ is a $\cal H$-$C^\ast$-algebra, then $A$
can also be considered to be a $\cal G$-$C^\ast$-algebra, and we have
a functorially induced morphism of $C^\ast$-categories $f_\star
A\rtimes_r {\cal G}\rightarrow A\rtimes_r {\cal H}$.

One can define the $K$-theory of $C^\ast$-categories; see
\cite{Mitch2.5} for details.  In particular, if $\cal A$ is a
$C^\ast$-category, there is an associated spectrum ${\mathbb K}({\cal
  A})$.  The assignment ${\cal A}\mapsto {\mathbb K}({\cal A})$ is
functorial.

We are now ready to state the main theorem of this article.

\begin{theorem} \label{2}
Let ${\mathbb E}^\%$ be a $G$-homotopy-invariant excisive functor from
the category of proper $G$-$CW$-complexes to the category of spectra.
Suppose we have a natural transformation $\alpha \colon {\mathbb E}^\%
(X)\rightarrow {\mathbb K}(A\rtimes_r \overline{X})$ such that the map
$$\alpha \colon {\mathbb E}^\% (G/H) \rightarrow {\mathbb
  K}(A\rtimes_r \overline{G/H})$$
is a stable equivalence for every finite subgroup, $H$, of the group
$G$.

Let $i\colon \overline{X}\rightarrow G$ be the obvious inclusion
functor.  Then the composite $i_\star \alpha \colon {\mathbb E}^\%
(X)\rightarrow {\mathbb K}(A\rtimes_r G)$ is the Baum-Connes assembly map.
\noproof
\end{theorem}

In order to prove the above theorem, we need to develop equivariant
$KK$-theory spectra of $C^\ast$-algebras in the world of groupoids.
This $KK$-theory must generalise equivariant $KK$-theory for groups,
and be related to crossed product $C^\ast$-categories.  The bulk of
this paper is devoted to the development of such a theory.

We should perhaps comment that Le Gall defines equivariant $KK$-theory
for groupoids in \cite{LG}.  However, Le Gall's approach is different to
ours, and it is not clear to the author how Le Gall's theory relates
to crossed product $C^\ast$-categories.  It is a potentially
interesting project to compare our theory with Le Gall's, but not a
project we explore in this paper.

\section{Preliminaries} \label{prelim}

Let $\mathbb F$ denote either the field of real numbers or the field
of complex numbers.  Recall that a {\em unital Banach category} over
the field $\mathbb F$ is a category, $\cal A$, in which every morphism
set $\Hom (A,B)_{\cal A}$ is a Banach space over the field $\mathbb F$, composition of morphisms
$$\Hom (B,C)_{\cal A}\times \Hom (A,B)_{\cal A}\rightarrow \Hom
(A,C)_{\cal A}$$
is bilinear, and the inequality
$$\| xy \| \leq \| x\| \|y \|$$
is satisfied for the norms of composable morphisms $x$ and $y$.

An {\em involution} on a Banach category $\cal A$ is a collection of
maps
$$\Hom (A,B)_{\cal A}\rightarrow \Hom (B,A)_{\cal A}$$
written $x\mapsto x^\star$
such that:

\begin{itemize}

\item $(\alpha x + \beta y)^\star = \overline{\alpha}x^\star +
  \overline{\beta}y^\star$ for all scalars $\alpha , \beta \in
  {\mathbb F}$ and morphisms $x,y\in \Hom (A,B)_{\cal A}$.

\item $(xy)^\star = y^\star x^\star$ for all composable morphisms $x$
  and $y$.

\item $(x^\star )^\star = x$ for every morphism $x$.

\end{itemize}

If $\cal A$ is a Banach category with involution, an invertible
morphism $u$ is called {\em unitary} if $u^{-1} = u^\star$.

The following definition comes from \cite{GLR} and \cite{Mitch2}

\begin{definition}
A unital Banach category with involution is called a {\em unital
  $C^\ast$-category} if for every morphism
  $x\in \Hom (A,B)_{\cal A}$, the product $x^\star x$ is a positive element of
  the Banach algebra $\Hom (A,A)_{\cal A}$, and the {\em $C^\ast$-identity}
$$\| x^\star x \| = \| x \|^2$$
holds.

A {\em non-unital $C^\ast$-category} is a collection of objects and
morphisms similar to a unital $C^\ast$-category except that there
need not exist identity morphisms $1\in \Hom (A,A)_{\cal A}$.
\end{definition}

We should perhaps comment that a non-unital $C^\ast$-category is not really a category, but rather an object with less structure which might be termed a {\em non-unital category}.

If $\cal A$ is a $C^\ast$-category, each
endomorphism set $\Hom (A,A)_{\cal A}$ is a $C^\ast$-algebra.  Conversely,
a $C^\ast$-algebra can be considered to be a $C^\ast$-category with
one object.

A {\em $C^\ast$-functor} between unital $C^\ast$-categories is a functor
$F\colon {\cal A}\rightarrow {\cal B}$ such that each map $F\colon \Hom
(A,B)_{\cal A}\rightarrow \Hom (F(A), F(B))_{\cal B}$ is linear, 
and $F(x^\star) = F(x)^\star$ for each morphism $x$ in the category
$\cal A$.  We similarly define $C^\ast$-functors between non-unital
$C^\ast$-categories.  It is proved in \cite{Mitch2} that any
$C^\ast$-functor is norm-decreasing, and therefore continuous, and if
faithful is an isometry.  Further, any $C^\ast$-functor has a closed image.

The category of small $C^\ast$-categories is formed by taking the
(non-unital) graded $C^\ast$-functors as
  morphisms.\footnote{A category is called {\em small} if the
  collection of objects is a set.  For set-theoretic reasons, one
  cannot form the category of all $C^\ast$-categories, whereas the
  category of small $C^\ast$-categories does make sense.}

\begin{example}
The category, ${\cal L}({\mathbb F})$, of all Hilbert spaces and
bounded linear operators is a $C^\ast$-category.  The involution is
defined by taking adjoints.
\end{example}

A $C^\ast$-functor $\rho \colon {\cal A}\rightarrow {\cal L}({\mathbb
  F})$ is termed a {\em representation} of the $C^\ast$-category $\cal
  A$.  It can be shown (see \cite{GLR, Mitch2}) that any small
$C^\ast$-category has a faithful, and therefore isometric
  representation.

For the applications we have in mind in this article it is necessary
to look at $C^\ast$-categories equipped with gradings.

\begin{definition}
A $C^\ast$-category $\cal A$ is said to be {\em graded} if we can
write each morphism set $\Hom (A,B)_{\cal A}$ as a direct sum
$$\Hom (A,B)_{\cal A} = \Hom (A,B)_0 \oplus \Hom (A,B)_1$$
of morphisms of degree $0$ and degree $1$ such that for composable morphisms $x$ and $y$ we have the formula
$$\deg (xy) = \deg (x) + \deg (y)$$

Here addition takes place modulo $2$.

A $C^\ast$-functor $F\colon {\cal A}\rightarrow {\cal B}$ between graded $C^\ast$-categories is termed a {\em graded $C^\ast$-functor} if
$$\deg (Fx) = \deg (x)$$
for every morphism $x$ in the category $\cal A$.
\end{definition}

As a special case of the above definition, we can speak of {\em
  graded $C^\ast$-algebras} and {\em morphisms} between graded
  $C^\ast$-algebras.  
The category of small graded $C^\ast$-categories
  is formed by taking the graded $C^\ast$-functors as
  morphisms.

We can consider an ungraded $C^\ast$-category to be equipped with the {\em trivial grading} defined by saying that every morphism is of degree $0$.  Our attitude is thus to view ungraded $C^\ast$-categories as special cases of graded $C^\ast$-categories.

There is a sensible notion of the spatial tensor product, ${\cal
  A}\gtp {\cal B}$, of graded $C^\ast$-categories $\cal A$ and $\cal
B$.  The objects are pairs, written $A\otimes B$, for objects $A\in
\Ob ({\cal A})$ and $B\in \Ob ({\cal B})$.  The morphism set $\Hom
(A\otimes B,A'\otimes B')_{{\cal A}\gtp {\cal B}}$ is a completion of
  the algebraic graded tensor product $\Hom (A,A')_{\cal A}\
  \hat{\odot}\ \Hom (B,B')_{\cal B}$.  See section 7 of \cite{Mitch2} and definition 2.7 of \cite{Mitch2.5} for details.

The main construction in \cite{Mitch2.5} is a functor, $\mathbb K$,
from the category of small graded $C^\ast$-categories to the category
of symmetric $\Omega$-spectra.  The spectrum ${\mathbb K}({\cal A})$ is
called the {\em $K$-theory spectrum} associated to the graded
$C^\ast$-category $\cal A$.\footnote{For an alternative construction
  of the $K$-theory spectrum of a $C^\ast$-category, see \cite{Jo2}.}
We define the $K$-theory group $K_n ({\cal A})$ to be the stable
homotopy group $\pi_n {\mathbb K}({\cal A})$.  If $A$ is a graded
$C^\ast$-algebra, the stable homotopy group, we recover from this
definition the $K$-theory groups $K_n (A)$ defined in \cite{vD1,
  vD2}.  In particular, when the $C^\ast$-algebra $A$ is trivially
graded, we can obtain the usual definition of $C^\ast$-algebra
$K$-theory in this way.

The $K$-theory of $C^\ast$-categories has many properties in common
with the $K$-theory of $C^\ast$-algebras.  A number of such elementary
properties are proved in the article \cite{Mitch2.5} including a
version of the Bott periodicity theorem involving Clifford
algebras.\footnote{See also \cite{ABS, Kar1, Kar2, Wo, vD2} for further details
  on this approach to the Bott periodicity theorem, at least for the
  $K$-theory of $C^\ast$-algebras.}

\begin{definition}
Let $p$ and $q$ be natural numbers.  Then we define the {\em $(p,q)$-Clifford algebra}, $\mathbb F_{p,q}$, to be the algebra over the field $\mathbb F$ generated by elements
$$\{ e_1 , \ldots , e_p , f_1 , \ldots , f_q \}$$
that pairwise anti-commute and satisfy the formulae
$$e_i^2 = 1 \qquad f_j^2 = -1$$ 
\end{definition}

The Clifford algebra ${\mathbb F}_{p,q}$ is a graded $C^\ast$-algebra;
the generators themselves are defined to be of degree $1$.  

\begin{theorem} \label{Bott}
Let $\cal A$ be a small graded $C^\ast$-category.  Then there is a natural stable equivalence of spectra
$$\Omega^q {\mathbb K}({\cal A}) \simeq \Omega^p {\mathbb K}({\cal A}\gtp \mathbb F_{p,q})$$
\noproof
\end{theorem}

Let $\{ {\cal A}_\lambda \ |\ \lambda \in \Lambda \}$ be a set
of small graded $C^\ast$-categories.  Then we can form the product, 
$\prod_{\lambda \in \Lambda} {\cal A}_\lambda$.  The objects are
collections of objects $\{ A_\lambda \in \lambda \in \Ob ({\cal
  A}_\Lambda ) \}$.  The morphism set $\Hom ( \{ A_\lambda \} , \{
B_\lambda \} )$ consists of all sets of morphisms $\{ x_\lambda \in
\Hom (A_\lambda , B_\lambda ) \ |\ \lambda \in \Lambda \}$ such that
the supremum $\sup \{ \| x_\lambda \| \ |\ \lambda \in \Lambda \}$ is
finite.

The following result is obvious from the construction of the
$K$-theory spectrum in \cite{Mitch2}.

\begin{proposition} \label{productcat}
Let $\{ {\cal A}_\lambda \ |\ \lambda \in \Lambda \}$ be a set of
small $C^\ast$-categories.  Define $\cal A$ to be the
$C^\ast$-category in which the set of objects is the union
$\bigcup_{\lambda \in \Lambda} \Ob (A_\lambda )$ and the morphism sets
are:
$$\Hom (A,B)_{\cal A} = \left\{ \begin{array}{ll}
\Hom (A,B)_{A_\lambda} & A,B \in \Ob ({\cal A}_\lambda ) \\
\{ 0 \} & A\in \Ob ({\cal A}_\lambda ), B\in \Ob ({\cal A}_\mu ), \
\lambda \neq \mu \\
\end{array} \right.$$

Then the $K$-theory spectra ${\mathbb K}({\cal A})$ and ${\mathbb
  K}(\prod_{\lambda \in \Lambda} {\cal A}_\lambda )$ are naturally
  stably equivalent.
\noproof
\end{proposition}

The other main property of $K$-theory that we need in this article is
a form of stability involving the objects of a $C^\ast$-category.

\begin{definition} \label{natiso}
Let $F, G\colon {\cal A}\rightarrow {\cal B}$ be graded $C^\ast$-functors
between unital graded $C^\ast$-categories.  Then a {\em natural isomorphism}
between $F$ and $G$ consists of a degree $0$ unitary morphism $U_A \in \Hom
(F(A),G(A))_{\cal B}$ for each object $A\in \Ob ({\cal A})$ such that
for every morphism $x\in \Hom (A,B)_{\cal A}$ the composites $U_B
F(x)$ and $F(x)U_A$ are equal.

A graded $C^\ast$-functor $F\colon {\cal A}\rightarrow {\cal B}$ between
unital $C^\ast$-categories is said to be an {\em equivalence} of
graded $C^\ast$-categories if there is a graded $C^\ast$-functor $G\colon {\cal B}\rightarrow {\cal A}$ such that the composites $FG$ and $GF$ are
naturally isomorphic to the identities $1_{\cal B}$ and $1_{\cal A}$ respectively.
\end{definition}

\begin{proposition} \label{natiso2}
Let $F\colon {\cal A}\rightarrow {\cal B}$ be an equivalence of small
graded $C^\ast$-categories.  Then the induced map $F_\star \colon {\mathbb
  K}({\cal A})\rightarrow {\mathbb K}({\cal B})$ is a stable
equivalence of $K$-theory spectra.
\noproof
\end{proposition}

In particular, a small graded unital $C^\ast$-category that is equivalent to a
$C^\ast$-algebra has the same $K$-theory.

We end our survey of results on the $K$-theory of $C^\ast$-categories
by indicating one way to define elements of the
initial space, ${\mathbb K}({\cal A})_0$, of the $K$-theory spectrum
${\mathbb K}({\cal A})$.  

Recall that a {\em right $\cal A$-module} over a $C^\ast$-category
$\cal A$ is a linear contravariant functor $\cal E$ from the category
$\cal A$ to the category of vector spaces.  It is similarly possible to define {\em left $\cal A$-modules}.

We use the notation
$$\eta x = {\cal E}(x) (\eta )$$
to denote the action of a morphism $x\in \Hom (A,B)_{\cal A}$ on a
vector $\eta \in {\cal E}(A)$.

\begin{definition}
The right $\cal A$-module $\cal E$ is called a {\em Hilbert $\cal
  A$-module} if it is equipped with a collection of bilinear maps
  $\langle -,-\rangle \colon {\cal E}(B)\times {\cal E}(A)\rightarrow
  \Hom (A,B)_{\cal A}$ such that:

\begin{itemize}

\item For all vectors $\eta \in {\cal E}(B)$, $\xi , \zeta \in {\cal
    E}(C)$, and morphisms $x,y\in \Hom (A,C)_{\cal A}$ we have the formula
$$\langle \eta , \xi x + \zeta y \rangle = \langle \eta ,\xi \rangle x
+ \langle \eta , \zeta \rangle y$$

\item $\langle \eta , \xi \rangle^\star = \langle \xi , \eta \rangle$

\item For each vector $\eta \in {\cal E}(A)$, the product $\langle
  \eta ,\eta \rangle$ is a positive element of the $C^\ast$-algebra
  $\Hom (A,A)$, and is zero only when the vector $\eta$ is zero.

\item Each vector space ${\cal E}(A)$ is complete with respect to the
  norm:
$$\| x \| = \| \langle x,x \rangle \|^\frac{1}{2}$$

\end{itemize}

\end{definition}

The collection of maps $\langle -,- \rangle \colon {\cal E}(B)\otimes
{\cal E}(A) \rightarrow \Hom (A,B)_{\cal A}$ is called an {\em inner
  product}.  If $\cal E$ is a Hilbert module over a
$C^\ast$-category $\cal A$ then each vector space ${\cal E}(A)$ is a
Hilbert module over the $C^\ast$-algebra $\Hom (A,A)$.  

Consider an object $A\in \Ob ({\cal A})$.  Then we have an associated
Hilbert $\cal A$-module $\Hom (-,A)_{\cal A}$.  The space associated
to the object $C\in \Ob ({\cal A})$ is the morphism set $\Hom (C
,A)_{\cal A}$.  The action of the category $\cal A$ is defined by
composition of morphisms, and the inner product is defined by the
formula
$$\langle x , y\rangle = x^\star y$$

There is an obvious notion of the direct sum, ${\cal E} \oplus {\cal
  F}$, of Hilbert $\cal A$-modules $\cal E$ and $\cal F$.  We can also
  define the direct sum of countably many Hilbert $\cal A$-modules;
  see definition 3.4 of \cite{Mitch3}.

We refer to two Hilbert $\cal A$-modules, $\cal E$ as {\em isomorphic} if there
  is a natural isomorphism of functors $T\colon {\cal E}\rightarrow
  {\cal F}$ such that
$$\langle \eta , \xi \rangle = \langle T\eta , T\xi \rangle$$
for all vectors $\eta \in {\cal E}(B)$ and $\xi \in {\cal E}(A)$.

\begin{definition}
A Hilbert $\cal A$-module $\cal E$ is called {\em finitely generated
  and projective} if there is a Hilbert $\cal A$-module ${\cal E}'$
  such that the direct sum ${\cal E}\oplus {\cal E}'$ is isomorphic to
  the direct sum of finitely many Hilbert $\cal A$-modules of the form
  $\Hom (-,A)_{\cal A}$.
\end{definition}

The following result is proved in \cite{Mitch3}.

\begin{proposition} \label{fgp}
Let $\cal A$ be a trivially graded unital $C^\ast$-category.  Then a
finitely generated projective Hilbert $\cal A$-module $\cal E$ defines
a canonical element of the initial space of the $K$-theory spectrum:
$[{\cal E}]\in {\mathbb K}({\cal A})_0$
\noproof
\end{proposition}

\section{Groupoid Actions} \label{actions}

Recall that a {\em groupoid} is a category in which every
morphism is invertible.  A group can be viewed as a groupoid with one object.
Taking this point of view, if $G$ is a discrete group, a {\em
  $G$-$C^\ast$-algebra} is a functor from the group $G$, viewed as a
category, to the category of $C^\ast$-algebras.  This idea
prompts the following definition.

\begin{definition}
Let $\cal G$ be a discrete groupoid.  Then a {\em $\cal G$-$C^\ast$-algebra}
is a functor from the groupoid $\cal G$ to the category of $C^\ast$-algebras.
\end{definition}

A $\cal G$-$C^\ast$-algebra as defined in \cite{LG} is also a $\cal
G$-$C^\ast$-algebra in the sense of the above definition.

If $A$ is a $\cal G$-$C^\ast$-algebra, let us write $A_a$
to denote the $C^\ast$-algebra associated to an object $a\in \Ob
({\cal G})$.  Then for each morphism $g\in \Hom (a,b)_{\cal G}$ we
have a morphism of $C^\ast$-algebras $g\colon A_a\rightarrow A_b$.

A $\cal G$-$C^\ast$-algebra $A$ is termed {\em unital} if every
$C^\ast$-algebra $A_a$ is unital, and the induced morphisms $g\colon
A_a \rightarrow A_b$ from the groupoid $\cal G$ all preserve the
unit.  A $\cal G$-$C^\ast$-algebra $A$ is termed {\em graded} if every
$C^\ast$-algebra $A_a$ is graded, and the induced morphisms $g\colon
A_a \rightarrow A_b$ from the groupoid $\cal G$ all preserve the grading.

\begin{definition}
A {\em $\cal G$-equivariant map} between $\cal G$-$C^\ast$-algebras
$A$ and $B$ is a natural transformation from the functor $A$ to the
functor $B$.
\end{definition}

More generally, let $f\colon {\cal G}\rightarrow {\cal H}$ be a
functor between groupoids, let $A$ be a $\cal G$-$C^\ast$-algebra,
and let $B$ be an $\cal H$-$C^\ast$-algebra.  Then an {\em equivariant
  map} $F\colon A\rightarrow B$ that {\em covers} the functor $f$ is a 
collection of morphisms of
$C^\ast$-algebras $F_a \colon A_a\rightarrow B_{f(a)}$ such that
$$F_b (gx) = f(g) F_a (x)$$
for every element $x\in A_a$ and morphism $g\in \Hom (a,b)_{\cal G}$.

If $A$ is a graded $\cal G$-$C^\ast$-algebra, and $B$ is a graded
$\cal H$-$C^\ast$-algebra, we insist that an equivariant map $F\colon
A\rightarrow B$ respects the gradings that are present.

\begin{example}
Let $f\colon {\cal G}\rightarrow {\cal H}$ be a functor between
groupoids, and let $A$ be an $\cal H$-$C^\ast$-algebra.  Then the
$\cal H$-$C^\ast$-algebra $A$ can also be considered a $\cal
G$-$C^\ast$-algebra; we associate the $C^\ast$-algebra $A_{f(a)}$ to
the object $a\in \Ob ({\cal G})$, and the morphism of
$C^\ast$-algebras $f(g)\colon
A_{f(a)}\rightarrow A_{f(b)}$ to the morphism $g\in \Hom (a,b)_{\cal
  G}$.

The collection of identity maps $1_a \colon A_{f(a)}\rightarrow
A_{f(a)}$ is an equivariant map that covers the functor $f$.
\end{example}

The above example will be important to us later on.

\begin{definition}
Let $A$ be a $\cal G$-$C^\ast$-algebra.  Then the {\em convolution
  category}, $A{\cal G}$, is the category with the same objects  as
  the groupoid $\cal G$ in which the morphism set $\Hom (a,b)_{{\cal
  G}A}$ consists of all formal sums:
$$x_1 g_1 + \cdots + x_n g_n$$
where $x_i \in A_b$ and $g_i \in \Hom (a,b)_{\cal G}$
\end{definition}

Composition of morphisms in the category ${\cal G}A$ is defined by the
formula
$$\left( \sum_i x_i g_i \right) \left( \sum_j y_j h_j \right) =
\sum_{i,j} x_i g_i (y_j )g_i h_j$$

Further, we have an involution
$$\left( \sum_i x_i g_i \right)^\star = \sum_i
  g_i^{-1}(x_i^\star ) g_i^{-1}$$

Note that the convolution category ${\cal G}A$ is non-unital unless
the $\cal G$-$C^\ast$-algebra $A$ is unital.

\begin{example} \label{trivial}
For a discrete groupoid $\cal G$, the {\em trivial $\cal
  G$-$C^\ast$-algebra} is defined by associating the scalar field
  $\mathbb F$ to each object $a\in \Ob ({\cal G})$, and the identity
  $1\colon {\mathbb F}\rightarrow {\mathbb F}$ to each morphism $g\in
  \Hom (a,b)_{\cal G}$.

The morphism set $\Hom (a,b)_{{\mathbb F}{\cal G}}$ in the convolution
$C^\ast$-category ${\mathbb F}{\cal G}$ consists of formal sums
$$\lambda_1 g_1 + \cdots + \lambda_n g_n$$
where $\lambda_i \in {\mathbb F}$ and $g_i \in \Hom (a,b)_{\cal G}$.

Composition and involution are defined by the formulae
$$\left( \sum_i \lambda_i g_i \right) \left( \sum_j \mu_j h_j \right) =
\sum_{i,j} \lambda_i \mu_j g_i h_j$$
and
$$\left( \sum_i \lambda_i g_i \right)^\star = \left( \sum_i
  \overline{\lambda_i} g_i^{-1} \right)$$
respectively.
\end{example}

Recall that we define ${\cal L}({\mathbb F})$ to be the
$C^\ast$-category of all Hilbert spaces and bounded linear operators
over the field $\mathbb F$.

\begin{definition}
A {\em unitary representation} of a groupoid $\cal G$ is a functor
$\rho \colon {\cal G}\rightarrow {\cal L}({\mathbb F})$ such that
$\rho (g^{-1}) = \rho (g)^\star$ for every morphism $g\in \Hom
(a,b)_{\cal G}$.
\end{definition}

We write $H_a$ to denote the Hilbert space associated to an object
$a\in \Ob ({\cal A})$.  The $C^\ast$-algebra of bounded linear operators
$T\colon H_a \rightarrow H_a$ is denoted ${\cal L}(H_a)$.

\begin{definition}
A {\em covariant representation} of a $\cal G$-$C^\ast$-algebra $A$ is a pair
$(\rho , \pi )$ consisting of a unitary representation $\rho \colon
{\cal G}\rightarrow {\cal L}({\mathbb F})$ 
together with representations $\pi \colon A_a \rightarrow {\cal
  L}(H_a)$ such that
$$\rho (g)\pi (x) = \pi (gx) \rho (g)$$
for every element $x\in A_a$ and morphism $g\in \Hom (a,b)_{\cal G}$.
\end{definition}

\begin{example} \label{regular}
Let $A$ be a $\cal G$-$C^\ast$-algebra.  Fix an object $a\in \Ob
({\cal G})$ and let $\alpha \colon A_a \rightarrow {\cal L}(H)$ be a
representation of the $C^\ast$-algebra $A_a$ on a Hilbert space $H$.
For each object $b\in \Ob ({\cal G})$, let $l^2 (a,b)$ be the Hilbert
space consisting of sequences
$(\eta_g )_{g\in \Hom (a,b)_{\cal G}}$ in the Hilbert space $H$ such
that the series $\sum_{g\in \Hom (a,b)_{\cal G}} \| \eta_g \|^2$
converges.

A groupoid element $h\in \Hom (b,c)_{\cal G}$ defines a unitary
operator $\rho (h)\colon l^2 (a,b)\rightarrow l^2 (a,c)$ by the formula
$$\rho (h)((\eta_g )_{g\in \Hom (a,b)_{\cal G}}) = ( \eta_{h^{-1}k})_{k\in
  \Hom (b,c)_{\cal G}}$$

We thus have a unitary representation, $\rho$, of the groupoid $\cal
G$ defined by mapping the object $b\in \Ob ({\cal G})$ to the Hilbert
space $l^2 (a,b)$, and the morphism $h\in \Hom (b,c)_{\cal G}$ to the above operator
$\rho (h)\colon l^2 (a,b)\rightarrow l^2 (b,c)$.

There are corresponding representations of the $C^\ast$-algebras $A_b$
defined by writing
$$\pi (x) ((\eta_g )_{g\in \Hom (a,b)_{\cal G}}) = (\alpha (g^{-1}(x))
\eta_g )_{g\in \Hom (a,b)_{\cal G}}$$

It is easy to verify that the formula
$$\rho (g)\pi (x) = \pi_b (g(x))\rho (g)$$
holds.  Therefore the pair $(\rho ,\pi )$ is a covariant
representation of the $\cal G$-$C^\ast$-algebra $A$.
\end{example}

A covariant representation of the type constructed in the above example is
called {\em regular}.

Associated to a covariant representation $(\rho ,\pi )$ we have a linear functor
$(\rho , \pi )_\star \colon A{\cal G}\rightarrow {\cal L}({\mathbb F})$
defined by mapping the object $a$ to the Hilbert space $H_a$ and the
morphism
$$x_1 g_1 + \cdots + x_n g_n \in \Hom (a,b)_{A{\cal G}}$$
to the bounded linear map
$$\pi (x_1) \rho (g_1) + \cdots + \pi (x_n) \rho(g_n) \colon H_a \rightarrow H_b$$

For any morphism $f\in \Hom (a,b)_{{\cal G}A}$, the
formula $(\rho , \pi )_\star (f^\star ) = (\rho , \pi )_\star (f)^\star$ holds.  We
express this formula by saying that the functor $(\rho , \pi )_\star$ {\em
  respects the involution.}

\begin{proposition} \label{allr}
Let $A$ be a unital $\cal G$-$C^\ast$-algebra.  Then every linear
functor $\alpha \colon A{\cal G}\rightarrow {\cal L}({\mathbb F})$
that respects the involution takes the form $(\rho , \pi )_\star$ for some
covariant representation $(\rho ,\pi )$.
\end{proposition}

\begin{proof}
Let $\alpha \colon A{\cal G}\rightarrow {\cal L}({\mathbb F})$ be a
linear functor that respects the involution.  Write $H_a = \alpha (a)$
for each object $a\in \Ob ({\cal G})$.  Then for any morphism $g\in
\Hom (a,b)$ we can define a unitary operator $\rho (g) \colon H_a
\rightarrow H_b$ by the formula
$$\rho (g) = \alpha (1_{A_b} g)$$
where $1_{A_b}$ is the identity element of the $C^\ast$-algebra
$A_b$.

We have a representation $\pi \colon A_a \rightarrow {\cal L}(H_a )$
defined by the formula
$$\pi (x) = \alpha (x1_a )$$
where $1_a \in \Hom (a,a)_{\cal G}$ is the identity morphism.  It is
easy to verify the formula
$$\rho(g)\pi (x) = \pi (g(x))\rho (g)$$
and the fact that $(\rho , \pi )_\star = \alpha$.
\end{proof}

The above result is also true in the non-unital case; we can prove it by
using approximate units for each $C^\ast$-algebra $A_a$.  However, we
do not need the non-unital result in this article, and therefore omit
the proof.

\begin{proposition}
Let $A$ be a $\cal G$-$C^\ast$-algebra.  Then we can define a norm on
the morphism sets of the convolution category ${\cal G}A$ by the
formula
$$\| \mu \|_\mathrm{max} = \sup \{ \| (\rho , \pi )_\star (\mu ) \| \
\textrm{$(\rho ,\pi )$ is a representation of $A$} \}$$
\end{proposition}

\begin{proof}
Consider a morphism
$$f = x_1 g_1 + \cdots + x_n g_n \in \Hom (a,b)_{{A\cal G}}$$

Observe that, for any representation $\pi$:
$$\| (\rho , \pi )_\star (f) \| \leq \| x_1 \| + \cdots + \| x_n \|$$

Hence the quantity $\| f \|_\mathrm{max}$ must be finite.  It is now easy to see that the function $f\mapsto \| f\|_\mathrm{max}$ is a norm so we are done.
\end{proof}

\begin{definition}
The {\em crossed product}, $A\rtimes {\cal G}$, is the Banach category
obtained by completing the morphism sets of the convolution category
$A{\cal G}$ with respect to the norm $\| - \|_\mathrm{max}$.
\end{definition}

The category $A\rtimes {\cal G}$ is equipped with an involution
inherited from the convolution category $A{\cal G}$.  It is
straightforward to verify the following result.

\begin{proposition}
The category $A\rtimes {\cal G}$ is a $C^\ast$-category.
\noproof
\end{proposition}

If $A$ is a graded $\cal G$-$C^\ast$-algebra, the crossed product
$A\rtimes {\cal G}$ can be graded by saying that a morphism
$$\sum_i x_i g_i \in \Hom (a,b)_{A\rtimes {\cal G}}$$
has degree $k$ if the elements $x_i \in A_b$ all have degree $k$.

There is another type of crossed product we need to consider.
To define it, observe that we can define a norm on the morphism sets of the
convolution category $A{\cal G}$ by the formula
$$\| \mu \|_r = \sup \{ \| (\rho , \pi )_\star (\mu ) \| \ \textrm{$(\rho ,\pi )$ is a regular representation of $A$} \}$$

\begin{definition}
The {\em reduced crossed product}, $A\rtimes_r {\cal G}$, is the
Banach category obtained by completing the morphism sets of the
convolution category $A{\cal G}$ with respect to the norm $\| - \|_r$.
\end{definition}

\begin{proposition}
The category $A\rtimes_r {\cal G}$ is a $C^\ast$-category.
\noproof
\end{proposition}

The reduced crossed product $A\rtimes_r {\cal G}$ is a graded
$C^\ast$-category when $A$ is a graded $\cal G$-$C^\ast$-algebra.

If $G$ is a discrete group we recover from the above definitions the
usual crossed product $C^\ast$-algebras $A\rtimes G$ and $A\rtimes_r G$.

Recall from example \ref{trivial} that for any groupoid $\cal G$ the {\em trivial $\cal G$-$C^\ast$-algebra} is defined by associating the scalar field
  $\mathbb F$ to each object $a\in \Ob ({\cal G})$, and the identity
  $1\colon {\mathbb F}\rightarrow {\mathbb F}$ to each morphism $g\in
  \Hom (a,b)_{\cal G}$.  

We have a corresponding convolution category
  ${\mathbb F}{\cal G}$, and crossed product $C^\ast$-categories
  ${\mathbb F}\rtimes {\cal G}$ and ${\mathbb F}\rtimes_r {\cal G}$.

\begin{definition}
The crossed product $C^\ast$-categories ${\mathbb F}\rtimes {\cal G}$
and ${\mathbb F}\rtimes_r {\cal G}$ are called the {\em maximal} and
{\em reduced} $C^\ast$-categories of the groupoid $\cal G$.  We
denote them by the symbols $C^\ast {\cal G}$ and $C^\ast_r {\cal G}$ respectively.
\end{definition}

The reduced and maximal $C^\ast$-categories of a groupoid were
originally defined without reference to crossed products in \cite{DL}
and \cite{Mitch2} respectively.  It is easy to see, using proposition
\ref{allr}, that the definitions given in these articles agree with
the above definition.

If $G$ is a group, we recover from the above definition
the usual maximal and reduced $C^\ast$-algebras, $C^\ast G$ and
$C^\ast_r G$, associated to the group $G$.

\begin{proposition} \label{crossfunctor}
Let $f\colon {\cal G}\rightarrow {\cal H}$ be a functor between
groupoids.  Let $A$ be a $\cal G$-$C^\ast$-algebra, $B$ be an $\cal
H$-$C^\ast$-algebra, and let $F\colon A\rightarrow B$ be an
equivariant map covering the functor $f$.

Then we have a functorially induced $C^\ast$-functor $F_\star \colon
A\rtimes {\cal G}\rightarrow B\rtimes {\cal H}$.
\end{proposition}

\begin{proof}
We begin by observing that we have an induced functor $F_\star \colon
A{\cal G}\rightarrow B{\cal H}$ between the convolution categories,
defined by writing $F_\star (a) =f(a)$ for each object $a\in \Ob
(A{\cal G})$ and
$$F_\star (x_1 g_1 + \cdots + x_n g_n ) = F(x_1 )f(g_1) + \cdots +
F(x_n)f(g_n)$$
for each morphism
$$x_1 g_1 + \cdots + x_n g_n \in \Hom (a,b)_{A{\cal G}}$$

Let $( \rho ,\pi )$ be a covariant representation
of the $\cal H$-$C^\ast$-algebra $B$.  Then
the pair $(\rho \circ f , \pi \circ F)$ is a representation of the
$\cal G$-$C^\ast$-algebra $A$.

For any morphism $\mu \in \Hom (a,b)_{{\cal G}A}$ we therefore have the
inequality
$$\| (\rho , \pi )_\star F_\star (\mu ) \| \leq \| \mu \|_\mathrm{max}$$

Hence $\| F_\star (\mu ) \|_\mathrm{max} \leq \| \mu \|_\mathrm{max}$ so the
functor $F_\star \colon A{\cal G}\rightarrow B{\cal G}$ is
continuous.  It therefore extends to a $C^\ast$-functor $F_\star
\colon A\rtimes {\cal G}\rightarrow B\rtimes {\cal H}$.

It is straightforward to check that the $C^\ast$-functor $F_\star
\colon A\rtimes {\cal G}\rightarrow B\rtimes {\cal G}$ depends
functorially on the equivariant map $F$, so we are done.
\end{proof}

If the equivariant map $F\colon A\rightarrow B$ respects the grading,
the induced $C^\ast$-functor $F_\star \colon A\rtimes {\cal
  G}\rightarrow B\rtimes {\cal H}$ is graded.

The $C^\ast$-categories $A\rtimes {\cal G}$ and $A\rtimes_r {\cal G}$
are not in general equal.  To see this, let $G$ be a group, and consider the trivial crossed products
  ${\mathbb F}\rtimes G = C^\ast G$ and ${\mathbb F}\rtimes_r G =
  C^\ast_r G$.  The group $C^\ast$-algebras $C^\ast G$ and $C^\ast_r
  G$ are equal if and only if the group $G$ is amenable.\footnote{See for example \cite{Pat} for details on amenability
  of groups.  In \cite{ADR} the idea of amenability for groupoids is
  analysed, which is of course relevant to the issue we are examining here.}

In fact, an example in \cite{DL} shows that it is impossible for the
assignment ${\cal G}\mapsto C^\ast_r {\cal G}$ to be functorial.
Thus the analogue of the above proposition for the reduced crossed
product is definitely false.  We do, however, have the following result.

\begin{proposition} \label{crossrfunctor}
Let $f\colon {\cal G}\rightarrow {\cal H}$ be a faithful functor between
groupoids.  Let $A$ be a $\cal G$-$C^\ast$-algebra, $B$ be an $\cal
H$-$C^\ast$-algebra, and let $F\colon A\rightarrow B$ be an injective
equivariant map covering the functor $f$.

Then we have a functorially induced $C^\ast$-functor $F_\star \colon
A\rtimes_r {\cal G}\rightarrow B\rtimes_r {\cal H}$.
\end{proposition}

\begin{proof}
As in proposition \ref{crossfunctor} we have an induced
functor $F_\star \colon A{\cal G}\rightarrow B{\cal H}$ between the
convolution categories.  We need to show that this functor is
continuous, and so extends to a $C^\ast$-functor $F_\ast \colon
A\rtimes_r {\cal G}\rightarrow B\rtimes_r {\cal H}$.

Choose an object $a\in \Ob ({\cal G})$ and a representation $\alpha
\colon B_{f(a)}\rightarrow {\cal L}(H)$.  Composition with the map $F$
yields a representation $\alpha F \colon A_a \rightarrow {\cal
  L}(H)$.  

Let $(\rho_\alpha , \pi_\alpha )$ and $(\rho_{\alpha F},
\pi_{\alpha F})$ be the induced regular representations of the
convolution categories $B{\cal H}$ and $A{\cal G}$ respectively,
defined as in example \ref{regular}.  Write $(B{\cal H})_\alpha$ and
$(A{\cal G})_{F\alpha }$ to denote the images of these regular
representations.  Then there is a faithful $C^\ast$-functor $F_\star
\colon (A{\cal G})_{F\alpha }\rightarrow (B{\cal H})_\alpha$ defined in the
obvious way.  Since any faithful $C^\ast$-functor is isometric, we
have the identity
$$\| (\rho_\alpha ,\pi_\alpha )_\star F_\star (\mu ) \| = \| (\rho
_{\alpha F},\pi_{\alpha F})_\star (\mu ) \|$$
for any morphism $\mu$ in the convolution category $A{\cal G}$.  

The definition of the norm in a reduced crossed product now gives us
the inequality
$$\| F_\star (\mu ) \|_r \leq \| \mu \|_r$$
and we are done.
\end{proof}

\begin{proposition} \label{canp}
Let $A$ be a $\cal G$-$C^\ast$-algebra.  Then we have a canonical
surjective $C^\ast$-functor $p\colon A\rtimes {\cal G}\rightarrow
A\rtimes_r {\cal G}$.
\end{proposition}

\begin{proof}
Observe that for any morphism $\mu$
in the convolution category we have the inequality
$$\| \mu \|_r \leq \| \mu \|_\mathrm{max}$$

Hence the identity functor on the convolution category, $1\colon
A{\cal G}\rightarrow A{\cal G}$, extends to a $C^\ast$-functor
$p\colon A\rtimes {\cal G}\rightarrow A\rtimes_r {\cal G}$.

As we remarked in section \ref{prelim}, any $C^\ast$-functor has a closed
image.\footnote{See corollary 4.9 in \cite{Mitch2}.}  Hence the image $p [\Hom (a,b)_{A\rtimes {\cal G}}]$ of a
morphism set in the $C^\ast$-category $A\rtimes {\cal G}$ is a
closed subset of the morphism set $\Hom (a,b)_{A\rtimes_r {\cal G}}$.
However, the image $p [\Hom (a,b)_{A\rtimes {\cal G}}]$ contains the
set $\Hom (a,b)_{A{\cal G}}$, which is a dense subset of the space
$\Hom (a,b)_{A\rtimes_r {\cal G}}$.  Therefore the $C^\ast$-functor
$p$ is surjective.
\end{proof}

The $C^\ast$-functor $p\colon A\rtimes {\cal G}\rightarrow A\rtimes_r
{\cal G}$ is natural in the category of faithful functors between
groupoids and injective equivariant maps.

\section{Equivariant $KK$-theory} \label{eqKK}

We begin our discussion of equivariant $KK$-theory by looking at
equivariant Hilbert modules.

\begin{definition}
Let $B$ be a $\cal G$-$C^\ast$-algebra.  Then a {\em $\cal G$-equivariant
Hilbert $B$-module} is a functor, $\cal E$, from the groupoid $\cal G$ to
the category of Banach spaces and invertible bounded linear maps such that:

\begin{itemize}

\item The space ${\cal E}_a$ associated to the object $a\in \Ob ({\cal G})$
  is a Hilbert $B_a$-module.

\item For every morphism $g\in \Hom (a,b)_{\cal G}$ we have the formula
$$g (\eta x) = (g\eta )(gx)$$
for all elements $\eta \in {\cal E}_a$ and $x\in B_a$.

\end{itemize}

\end{definition}

A $\cal G$-equivariant Hilbert $B$-module $\cal E$ is said to be {\em countably
  generated} if each Hilbert $B_a$-module ${\cal E}_a$ is countably
  generated.

If $B$ is a graded $\cal G$-$C^\ast$-algebra, a $\cal G$-equivariant
Hilbert $B$-module $\cal E$ is referred to as {\em graded} if each Hilbert
$B_a$-module ${\cal E}_a$ is graded, and the maps $g\colon {\cal
  E}_a\rightarrow {\cal E}_b$ coming from the groupoid $\cal G$ all
respect the grading.

\begin{example} \label{Aeq}
Let $B$ be a $\cal G$-$C^\ast$-algebra.  Then $B$ itself can be
considered to be a $\cal G$-equivariant Hilbert module.  We have inner
products $\langle -,-\rangle \colon B_a\times B_a \rightarrow B_a$
defined by the formula
$$\langle x,y\rangle = x^\star y$$
\end{example}

If $B$ is a graded $\cal G$-$C^\ast$-algebra, then $B$ is also graded
as a $\cal G$-equivariant Hilbert $B$-module.

\begin{definition}
Let $\cal E$ be a $\cal G$-equivariant Hilbert $B$-module.  We write ${\cal L}_{\cal G}({\cal E})$ to denote the $\cal G$-$C^\ast$-algebra
which associates the $C^\ast$-algebra ${\cal L}({\cal E}_a )$ to the object
$a\in \Ob ({\cal G})$.  The $\cal G$-action is defined by the formula
$$g(T\eta ) = (gT)(g\eta)$$
\end{definition}

The $\cal G$-$C^\ast$-algebra ${\cal L}_{\cal G}({\cal E})$ is graded in the
obvious way when $\cal E$ is a $\cal G$-equivariant graded Hilbert $B$-module.

For graded $C^\ast$-algebras $A$ and $B$ we define a {\em graded Hilbert
  $(A,B)$-bimodule} to be a countably generated graded Hilbert $B$-module
  $\cal F$ equipped with a morphism $\phi \colon A\rightarrow {\cal L}({\cal F})$.  The right
  $B$-module $\cal F$ is thus also a left $A$-module, with $A$-action:
$$x \eta = \phi (x) (\eta )$$

We will usually drop explicit mention of the morphism $\phi$ from our
notation.  Let $\cal E$ be a graded Hilbert $A$-module.  Then we can form the {\em algebraic tensor product} ${\cal E}\odot_A {\cal F}$; it is the right $B$-module
generated by elementary tensors $\eta \otimes \xi$, where $\eta \in
{\cal E}$
and $\xi \in {\cal F}$, subject to the relation
$$\eta x\otimes \xi = \eta \otimes x\xi$$
for all elements $x\in A$.  The {\em inner tensor product},
${\cal E}\otimes_A {\cal F}$, is the completion of the algebraic tensor product
${\cal E}\odot_A {\cal F}$ with respect to the norm defined by the inner product
$$\langle \eta \otimes \xi , \eta' \otimes \xi' \rangle = \langle \xi
,\langle \eta , \eta' \rangle \xi' \rangle$$

See \cite{Lan} for further details.

\begin{definition}
Let $A$ and $B$ be $\cal G$-$C^\ast$-algebras.  Then a {\em $\cal
  G$-equivariant graded Hilbert $(A,B)$-bimodule} is a $\cal
  G$-equivariant graded Hilbert $B$-module, $\cal F$, equipped with a $\cal G$-equivariant map $\phi \colon A\rightarrow {\cal L}_{\cal G}({\cal F})$.
\end{definition}

As in the non-equivariant case, we usually drop explicit mention of the
morphism $\phi$ from our notation.  Observe that if $\cal F$ is a
$\cal G$-equivariant graded Hilbert $(A,B)$-bimodule, each space
${\cal F}_a$ is a graded Hilbert $(A_a,B_a)$-bimodule.

We call a $\cal G$-equivariant graded Hilbert $(A,B)$-bimodule {\em
  countably generated} if it is countably generated as a $\cal
  G$-equivariant graded Hilbert $B$-module.

\begin{example} \label{onebiq}
Let $A$ and $B$ be graded $\cal G$-$C^\ast$-algebras and suppose we have a
$\cal G$-equivariant map $\phi \colon A\rightarrow B$.
According to example \ref{Aeq} the $\cal G$-$C^\ast$-algebra $B$ is
itself a graded $\cal G$-equivariant Hilbert $B$-module.  Given an element $x\in A_a$ we have an
operator $\phi (x) \colon B_a \rightarrow B_a$ defined by
multiplication.  The formula $g(\phi (x)y) = \phi (gx)(gy)$ is
satisfied.  The map $\phi$ can therefore by considered a $\cal
G$-equivariant map $\phi \colon A\rightarrow {\cal L}_{\cal G}(B)$,
and the $\cal G$-$C^\ast$-algebra $B$ is itself a $\cal
G$-equivariant Hilbert $(A,B)$-bimodule.
\end{example}

\begin{definition}
Let $\cal E$ be a $\cal G$-equivariant graded Hilbert $A$-module, and let
$\cal F$ be a $\cal G$-equivariant graded Hilbert $(A,B)$-bimodule.  Then the
{\em inner tensor product}, ${\cal E}\otimes_A {\cal F}$, is the $\cal
G$-equivariant graded Hilbert
$B$-module in which the module $({\cal E}\otimes_A {\cal F})_a$ is the tensor
product ${\cal E}_a \otimes_{A_a} {\cal F}_a$.  The $\cal G$-action is defined by
the formula
$$g (\eta \otimes \xi ) = g\eta \otimes g\xi$$
\end{definition}

The inner tensor product ${\cal E}\otimes_A {\cal F}$ is countably
generated if the $\cal G$-equivariant Hilbert $A$-module $\cal E$ and
$\cal G$-equivariant Hilbert $(A,B)$-bimodule $\cal F$ are countably generated.

There is another type of tensor product of Hilbert modules that we
will need in our calculations.  Recall that if $A$ and $B$ are
graded $C^\ast$-algebras, $\cal E$ is a graded Hilbert $A$-module, and
$\cal F$ is a graded Hilbert
$B$-module, we can define a Hilbert $A\gtp B$-module ${\cal
  E}\otimes {\cal F}$.
It is the completion of the tensor product of vector spaces ${\cal
  E}\odot {\cal F}$
with respect to the norm defined by the inner product
$$\langle \eta \otimes \xi , \eta' \otimes \xi' \rangle = \langle \eta
,\eta' \rangle \otimes \langle \xi ,\xi' \rangle$$

The relevant technical details can be found in \cite{Lan}.  We can define a grading by specifying the degree of elementary tensors:
$$\deg (\eta \otimes \xi ) = \deg (\eta ) + \deg (\xi )$$

Here addition takes place modulo $2$.

\begin{definition} \label{outer}
Let $A$ and $B$ be graded $\cal G$-$C^\ast$-algebras.  Then we define the
{\em tensor product}, $A\gtp B$, to be the $\cal
G$-$C^\ast$-algebra in which the $C^\ast$-algebra $(A\gtp B)_a$ is
equal to the tensor product $A_a\gtp B_a$.  The $\cal G$-action is
defined by the formula
$$g (x \otimes y ) = gx \otimes gy$$

If $\cal E$ is a $\cal G$-equivariant graded Hilbert $A$-module, and
$\cal F$ is a
$\cal G$-equivariant graded Hilbert $B$-module, we define the {\em outer tensor
  product}, ${\cal E}\otimes {\cal F}$, to be the $\cal G$-equivariant
graded Hilbert
$A\otimes B$-module where the module $({\cal E}\otimes {\cal F})_a$ is
the tensor product ${\cal E}_a \otimes {\cal F}_a$.  The $\cal
G$-action is defined by the formula
$$g (\eta \otimes \xi ) = g\eta \otimes g\xi$$
\end{definition}

The outer tensor product of two countably generated $\cal
G$-equivariant Hilbert modules is countably generated.

We need one final definition before we are ready to look at
equivariant $KK$-theory.

\begin{definition}
Let $A$ be a $\cal G$-$C^\ast$-algebra, and let $\cal E$ and ${\cal
  E}'$ be $\cal G$-equivariant Hilbert $A$-modules.  Then a {\em bounded operator} $T\colon {\cal E}\rightarrow {\cal E}'$
is a collection, $T$, of operators $T_a \colon {\cal
  E}_a\rightarrow {\cal E}'_a$ such that the norm
$$\| T \| = \sup \{ \| T_a \| \ |\ a\in \Ob ({\cal G}) \}$$
is finite.
\end{definition}

Note that we make no assumptions here concerning equivariance.  In the
graded case, we say a bounded operator has degree $k$ if each operator
$T_a \colon {\cal E}_a \rightarrow {\cal E}'_a$ has degree $k$.

If $A$ and $B$ are graded $\cal G$-$C^\ast$-algebras, and $\cal E$ and
${\cal E}'$ are $\cal G$-equivariant graded Hilbert $(A,B)$-bimodules,
we call a collection, $T$, of maps
$T_a\colon {\cal E}_a\rightarrow {\cal E}'_a$ a {\em bounded operator}
if it is a bounded operator between $\cal G$-equivariant Hilbert $B$-modules
in the sense of the above definition.

Let $T$ be a bounded operator between $\cal G$-equivariant graded
Hilbert $(A,B)$-bimodules, and let $x\in A_a$.  Then we define the {\em graded commutator}
$$[x,T] = xT_a - (-1)^{\deg (x) \deg (T)} T_a x$$

This formula only makes sense when the degree of the element $x$
and the operator $T$ are defined.  However, we can extend the 
definition of the graded commutator by requiring it to be linear in
each variable.

\begin{definition}
Let $A$ and $B$ be $\cal G$-$C^\ast$-algebras. Then a {\em $\cal
  G$-equivariant Kasparov $(A,B)$-cycle} is a pair $({\cal E},T)$,
  where $\cal E$ is a countably generated $\cal G$-equivariant Hilbert
  $(A,B)$-bimodule, and $T\colon {\cal E} \rightarrow {\cal E}$ is a bounded operator such that the operators
$$x(T_a-T_a^\star) \qquad x(T_a^2 - 1) \qquad [x , T] \qquad x(gT_b -
T_a g)$$
are compact for all elements $x\in A_a$ and morphisms $g\in \Hom
(b,a)_{\cal G}$.
\end{definition}

In the above definition, the operators defined by the various formulae
are just operators between (non-equivariant) Hilbert modules over
$C^\ast$-algebras.  We use the standard $C^\ast$-algebraic notion of
such an operator being compact.

An element of equivariant $KK$-theory is a certain equivalence class
of equivariant Kasparov cycles.  

\begin{definition}
Let $({\cal E},T)$ and $({\cal E}',T')$ be $\cal G$-equivariant Kasparov $(A,B)$-cycles.  Then the {\em direct sum} is the Kasparov cycle
$$({\cal E},T)\oplus ({\cal E}',T') = ({\cal E}\oplus {\cal E}',T\oplus T')$$
\end{definition}

\begin{definition}
\mbox{}

\begin{itemize}

\item[$\bullet$] A $\cal G$-equivariant Kasparov $(A,B)$-cycle $({\cal E},T)$ is called {\em degenerate} if the operators
$$x(T_a-T_a^\star) \qquad x(T_a^2 - 1) \qquad [x , T] \qquad x(gT_b -
T_a g)$$
are equal to zero for all elements $x\in A_a$ and morphisms $g\in \Hom
(b,a)_{\cal G}$.

\item[$\bullet$] An {\em operator homotopy} between $\cal
  G$-equivariant Kasparov $(A,B)$-cycles $({\cal E},T)$ and $({\cal E},T')$ is a norm-continuous path $({\cal E},T_t)$ of Kasparov cycles such that $T_0 = T$ and $T_1 = T'$.

\item[$\bullet$] Two $\cal G$-equivariant Kasparov $(A,B)$-cycles
  $({\cal E}_1,T_1)$ and $({\cal E}_2, T_2)$ are called {\em
    equivalent} if there are degenerate Kasparov cycles $({\cal E}_1',T_1')$ and $({\cal E}_2',T_2')$ such that the direct sums $({\cal E}_1,T_1)\oplus({\cal E}_1',T_1')$ and $({\cal E}_2, T_2)\oplus({\cal E}_2', T_2')$ are operator homotopic.

\end{itemize}

\end{definition}

We write $[({\cal E}, T)]$ to denote the equivalence class of a $\cal
G$-equivariant Kasparov $(A,B)$-cycle $({\cal E},T)$, and $KK_{\cal G}(A,B)$ to denote the set of equivalence classes.

\begin{proposition}
The set $KK_{\cal G}(A,B)$ is an Abelian group with an operation defined by taking the direct sum of Kasparov cycles.  
\end{proposition}

\begin{proof}
It is easy to check that the set $KK_{\cal G}(A,B)$ is an Abelian
semigroup, with identity element $[({\cal E},T)]$ where $({\cal E},T)$
is any degenerate $\cal G$-equivariant Kasparov $(A,B)$-cycle.

If ${\cal E}$ is a $\cal G$-equivariant Hilbert $B$-module, with
grading ${\cal E}_a= ({\cal E}_a)_0 \oplus ({\cal E}_a)_1$, define
${\cal E}^\mathrm{op}$ to be the $\cal G$-equivariant Hilbert
$B$-module with the opposite grading.  If we have a $\cal
G$-equivariant map $\phi \colon A\rightarrow {\cal L}({\cal E})$ we
can define a $\cal G$-equivariant map $\phi^\mathrm{op} \colon
A\rightarrow {\cal L}({\cal E}^\mathrm{op})$ by writing
$\phi^\mathrm{op} ( x_0 + x_1) = \phi (x_0 - x_1)$ for elements $x_0 ,
x_1 \in A_a$ of degrees $0$ and $1$ respectively.

Consider a $\cal G$-equivariant Kasparov $(A,B)$-cycle $({\cal
  E},T)$.  The pair $({\cal E}^\mathrm{op} , -T)$ is a $\cal
  G$-equivariant Kasparov $(A,B)$-cycle.
We can define an operator homotopy between the Kasparov cycle $({\cal
  E},T)\oplus ({\cal E}^\mathrm{op},-T)$ and the degenerate cycle $\left( {\cal E}\oplus \tilde{\cal E} , \left( \begin{array}{cc}
0 & 1 \\
1 & 0 \\
\end{array} \right) \right)$ by the formula
$$G_\theta = \left( \begin{array}{cc}
T\cos \theta & \sin \theta \\
\sin \theta & -T\cos \theta \\
\end{array} \right) \qquad \theta \in [0, \frac{\pi}{2}]$$

Hence
$$[({\cal E},T)] + [(\tilde{\cal E}, -T)] = 0$$
and we have proved that the set $KK_{\cal G}({\cal A},{\cal B})$ is an Abelian group.
\end{proof}

When $G$ a group, we recover from the above definition the usual
equivariant $KK$-theory groups, as defined by Kasparov in
\cite{Kas3}.  

There is a small technical point here that we should mention.  It is
usual to define the equivariant $KK$-theory groups by looking at an
equivalence relation called {\em homotopy} on Kasparov cycles, rather
than the equivalence relation we have used in our more general
definition.  However, these relations turn out to be the same; see
remark 5.11 (2) in \cite{BS}.

Let us call a $\cal G$-$C^\ast$-algebra $A$ {\em $\sigma$-unital} if
each $C^\ast$-algebra $A_a$ has a countable approximate unit.

\begin{proposition} \label{func0}
Let $A$ and $B$ be $\sigma$-unital $\cal G$-$C^\ast$-algebras.  The
group $KK_{\cal G}(A,B)$ is contravariantly functorial in the variable
$A$ and covariantly functorial in the variable $B$.
\end{proposition}

\begin{proof}
Let $({\cal E},T)$ be a $\cal G$-equivariant Kasparov $(A,B)$-cycle,
and let $F\colon B\rightarrow B'$ be a $\cal G$-equivariant map.  

Then the $\cal G$-$C^\ast$-algebra $B'$ is itself a countably
generated $\cal G$-equivariant graded Hilbert $(B,B')$-bimodule according to example
\ref{onebiq}.\footnote{We need the $\cal G$-$C^\ast$-algebra $B'$ to
  be $\sigma$-unital for the $\cal G$-equivariant Hilbert
  $(B,B')$-module $B'$ to be countably generated.}  
We can therefore form the inner tensor product ${\cal
  E}\otimes_B B'$.  This inner tensor product is a $\cal
G$-equivariant graded Hilbert $(A,B')$-bimodule since we have an $A$-action
defined by writing
$$x (\eta \otimes y) = (x\eta )\otimes y$$
where $x\in A_a$, $\eta \in {\cal E}_a$, and $y\in B'_a$.

Further, there is a bounded operator
$T\otimes 1 \colon {\cal E}\otimes_B B'\rightarrow {\cal E}'\otimes_B
B'$ given by the formula
$$(T\otimes 1)(\eta \otimes y) = (T\eta )\otimes y$$

It is easy to check that we have a functorially induced map $F_\star \colon
KK_{\cal G}(A,B)\rightarrow KK_{\cal G}(A,B')$ defined by writing
$$F_\star [({\cal E},T)] = [({\cal E}\otimes_B B',T\otimes 1)]$$

Now consider a $\cal G$-equivariant map $G\colon A'\rightarrow A$.
Suppose that the action of the $\cal G$-$C^\ast$-algebra $A$ is
defined on the Hilbert $B$-module $\cal E$ by the equivariant map
$\phi \colon A\rightarrow {\cal L}_{\cal G}({\cal E})$.  Then we can
form a $\cal G$-equivariant graded Hilbert $(A',B)$-bimodule $G^\star ({\cal E})$.  The module $G^\star ({\cal E})$ is equal to
the module $\cal E$ as a $\cal G$-equivariant graded Hilbert
$B$-module, and the $A'$-action is defined by the equivariant map $\phi \circ G
\colon A'\rightarrow {\cal L}_{\cal G}G^\star ({\cal E})$.

We have a functorially induced map $G^\star \colon
KK_{\cal G}(A,B)\rightarrow KK_{\cal G}(A',B)$ defined by the formula
$$G^\star [({\cal E},T)] = [(G^\star({\cal E}) , T)]$$
\end{proof}

In \cite{Hi2} the $K$-homology of a $C^\star$-algebra $A$ is defined
in terms of the ordinary $K$-theory of a `dual algebra' constructed
from $A$.  We can extend this approach to define the equivariant
$KK$-theory groups $KK^{-n}_{\cal G}({\cal A},{\cal B})$ for $\cal
G$-$C^\ast$-algebras $A$ and $B$ in terms of the ordinary $K$-theory
of some `dual $C^\ast$-category'.  

Our definitions and methods here are modeled on the approach to the
$KK$-theory of $C^\ast$-categories in \cite{Mitch3}.

\begin{definition}
We write ${\cal D}_{\cal G}(A,B)$ to denote the category of
countably generated $\cal G$-equivariant Hilbert $(A,B)$-bimodules and
bounded operators $T\colon {\cal E}\rightarrow {\cal E}'$ such that
the operators $[x, T]$ and $x(gT_a - T_b g)$ are compact for all
elements $x\in A_a$ and morphisms $g\in \Hom (a,b)_{\cal G}$.  We
write ${\cal K}{\cal D}_{\cal G}(A,B)$ to
denote the (non-unital) subcategory consisting of bounded operators
$T\colon {\cal E}\rightarrow {\cal E}'$ such that the composites
$xT_a$ and $T_ax$ are compact operators for all elements $x\in
A_a$.
\end{definition}

A {\em $C^\ast$-ideal}, $\cal J$, in a $C^\ast$-category $\cal A$ is a
$C^\ast$-subcategory such that the composite of a morphism in the
category $\cal J$ and a morphism in the category $\cal A$ belongs to
the category $\cal J$.  One can form the quotient, ${\cal A}/{\cal
  J}$, of a $C^\ast$-category by a $C^\ast$-ideal; see \cite{Mitch2}
for details.

A straightforward calculation tells us that the category 
${\cal D}_{\cal G}(A,B)$ is a graded $C^\ast$-category and
the subcategory ${\cal K}{\cal D}_{\cal G}(A,B)$ is a 
$C^\ast$-ideal.  We can therefore form the quotient
$${\cal Q}{\cal D}_{\cal G}(A,B) = {\cal D}_{\cal G}(A,B) / {\cal K}{\cal
  D}_{\cal G}(A,B)$$

The following result is proved in the same way as proposition \ref{func0}.

\begin{proposition} \label{func1}
Let $A$ and $B$ be $\sigma$-unital $\cal G$-$C^\ast$-algebras.  Then
the graded $C^\ast$-category ${\cal Q}{\cal D}_{\cal G}(A,B)$ is
contravariantly functorial in the variable $A$ and covariantly
functorial in the variable $B$.
\noproof
\end{proposition}

The following results are proved in exactly the same way as theorem
4.9 and lemma 4.10 in \cite{Mitch3}.  We do not repeat the work here.

\begin{theorem}
Let $A$ and $B$ be $\sigma$-unital $\cal G$-$C^\ast$-categories.  There is a natural isomorphism $K_1 {\cal QD}_{\cal G}(A,B)\cong KK(A,B)$.
\noproof
\end{theorem}

\begin{lemma}
Suppose that $A$ and $B$ be graded $\cal G$-$C^\ast$-algebras.
Let $p,q\in \mathbb N$.  Then there is a natural isomorphism
$$K_1{\cal QD}_{\cal G}(A,B\gtp \mathbb F_{p,q})\cong K_1({\cal
  QD}_{\cal G}(A,B)\gtp \mathbb F_{p,q})$$
\noproof
\end{lemma}

By the Bott periodicity theorem it therefore makes sense to define
further equivariant $KK$-theory groups by the formula
$$KK^{p-q}_{\cal G} (A,B) = KK_{\cal G}(A,B\gtp \mathbb F_{p,q})$$

We have natural isomorphisms
$$KK^{p-q}_{\cal G}(A,B)\cong K_1({\cal QD}_{\cal G}(A,B)\gtp \mathbb
F_{p,q})\cong K_{1-(p-q)}{\cal QD}_{\cal G}(A,B)$$

\begin{definition}
Let $\cal G$ be a discrete groupoid, and let $A$ and $B$ be
$\sigma$-unital $\cal G$-$C^\ast$-categories.  Then we define the 
{\em $\cal G$-equivariant $KK$-theory spectrum}
$${\mathbb K}{\mathbb K}_{\cal G}(A,B) = \Omega {\mathbb K} {\cal
  Q}{\cal D}_{\cal G}(A,B)$$
\end{definition}

According to proposition \ref{func1} the $KK$-theory spectrum ${\mathbb
  K}{\mathbb K}_{\cal G}(A,B)$ is contravariantly functorial in
  the variable $A$ and covariantly functorial in the variable
  $B$.  There is another type of functoriality that we need to
  consider, depending this time on the groupoid $\cal G$.

Let $f\colon {\cal H}\rightarrow {\cal G}$ be a functor between
groupoids, and let $A$ be a graded $\cal G$-$C^\ast$-algebras.  Abusing
notation, $A$ can also be considered to be a graded $\cal
H$-$C^\ast$-algebra.  We associate the graded $C^\ast$-algebra $A_{f(a)}$ to
the object $a\in \Ob ({\cal G})$ and the morphism $f(g) \colon
A_{f(a)}\rightarrow A_{f(a)}$ to the element $g\in \Hom (a,b)_{\cal
  H}$.

\begin{proposition} \label{restriction}
There is a functorially induced graded $C^\ast$-functor $f^\star \colon
{\cal Q}{\cal D}_{\cal G} (A,B)\rightarrow {\cal Q}{\cal D}_{\cal H}
(A,B)$.  The $C^\ast$-functor $f^\star$ is natural in the variables $A$ and
$B$.
\end{proposition}

\begin{proof}
Consider a $\cal G$-equivariant graded Hilbert $(A,B)$-bimodule $\cal
E$.  Then we have an $\cal H$-equivariant graded Hilbert $B$-bimodule $f^\star (\cal E )$
defined by associating the Hilbert $B_{f(a)}$-module ${\cal E}_{f(a)}$
to the object $a\in \Ob ({\cal H})$.  The action of the groupoid $\cal
H$ is defined by the formula $g\eta = f(g) \eta$.  The action of the
$\cal H$-$C^\ast$-algebra $A$ is the same as the action of the $\cal
G$-$C^\ast$-algebra $A$ on the original bimodule $\cal E$.

Let $\cal E$ and ${\cal E}'$ be graded Hilbert $(A,B)$-modules, and let
$T\colon {\cal E}\rightarrow {\cal E}'$ be a bounded operator.  Then
we have a bounded operator $f^\star (T) \colon f^\star ({\cal E})
\rightarrow f^\star ({\cal E}' )$ defined by the formula $f^\star
(T)(\eta ) = T\eta$.

If the operator $T$ is a morphism in the category ${\cal D}_{\cal
  G}(A,B)$, the operator $f^\star (T)$ is a morphism in the category
  ${\cal D}_{\cal H}(A,B)$.  If the operator $T$ is a morphism in the
  category ${\cal K}{\cal D}_{\cal G}(A,B)$, the operator $f^\star (T)$ is a morphism in the category
  ${\cal K}{\cal D}_{\cal H}(A,B)$.  We therefore
  have a functorially induced graded $C^\ast$-functor $f^\star \colon {\cal
  Q}{\cal D}_{\cal G}(A,B)\rightarrow {\cal Q}{\cal D}_{\cal
  H}(A,B)$.  Naturality of this induced $C^\ast$-functor in the
  variables $A$ and $B$ is easy to check.
\end{proof}

In particular, at the level of $K$-theory, we have a map
$$f^\star \colon {\mathbb K}{\mathbb K}_{\cal G}(A,B)\rightarrow
{\mathbb K}{\mathbb K}_{\cal H}(A,B)$$

\begin{definition} 
The induced map
$$f^\star \colon {\mathbb K}{\mathbb K}_{\cal G}(A,B)\rightarrow
{\mathbb K}{\mathbb K}_{\cal H}(A,B)$$
is called the {\em restriction map}.
\end{definition}

When $f$ is a group homomorphism, we recover from the above definition
the usual restriction maps in equivariant $KK$-theory.

\begin{proposition} \label{gequivalence}
Let $f\colon {\cal H}\rightarrow {\cal G}$ be an equivalence of
groupoids.  Let $A$ and $B$ be unital $\cal G$-$C^\ast$-algebras.  
Then the restriction map
$$f^\star \colon {\mathbb K}{\mathbb K}_{\cal G}(A,B)\rightarrow
{\mathbb K}{\mathbb K}_{\cal H}(A,B)$$
is a stable equivalence of spectra.
\end{proposition}

\begin{proof}
By proposition \ref{natiso2} it suffices to show that the $C^\ast$-functor
$$f^\star \colon {\cal Q}{\cal D}_{\cal G}(A,B)\rightarrow {\cal
  Q}{\cal D}_{\cal H}(A,B)$$
is an equivalence of $C^\ast$-categories.  Let $g\colon {\cal
  G}\rightarrow {\cal H}$ be a functor such that we have natural
isomorphisms $F\colon fg\rightarrow 1_{\cal G}$ and $G\colon
gf\rightarrow 1_{\cal H}$ respectively.

Consider an object $a\in \Ob ({\cal G})$.  Then we have a morphism
$F_a \in \Hom (fg(a),a)_{\cal G}$ determined by the natural isomorphism $F$.
If $\cal E$ is a $\cal G$-equivariant Hilbert $B$-module, the action
of the groupoid $\cal G$ gives us a unitary operator
$$F_a \colon g^\star f^\star ({\cal E})\rightarrow {\cal E}$$

It is easy to check that we have an induced natural isomorphism
$$F_\star \colon g^\star f^\star \rightarrow 1_{{\cal QD}_{\cal
    G}(A,B)}$$

Similarly we obtain a natural isomorphism
$$G_\star \colon f^\star g^\star \rightarrow 1_{{\cal QD}_{\cal
    H}(A,B)}$$

Therefore the $C^\ast$-functor $f^\star \colon {\cal Q}{\cal D}_{\cal G}(A,B)\rightarrow {\cal
  Q}{\cal D}_{\cal H}(A,B)$
is an equivalence of $C^\ast$-categories and we are done.
\end{proof}

\section{Descent}

There is a canonical {\em descent map} from the equivariant
$KK$-theory of $\cal G$-$C^\ast$-algebras to the $KK$-theory of the
associated crossed product $C^\ast$-categories.  Before we construct
it, we need to review the $KK$-theory of $C^\ast$-categories, as
defined in \cite{Mitch3}.

The $KK$-theory of $C^\ast$-categories is constructed by considering certain operators between countably generated Hilbert modules.

We begin by making the notion of countably generated precise in this context.
Let $\cal B$ be a $C^\ast$-category.  Recall that a right $\cal B$-module, $\cal E$, is said to be {\em countably
  generated} if there is a countable set
$$\Omega \subseteq \bigcup_{A\in \Ob ({\cal B})} {\cal E}(A)$$
such that for each object $A\in \Ob ({\cal B})$, every element of the vector
space ${\cal E}(A)$ is a finite linear combination of elements of the
form $\eta x$, where $x\in \Hom (A,B)_{\cal B}$ and $\eta \in \Omega \cap {\cal E}(B)$.

\begin{definition}
A Hilbert $\cal B$-module $\cal E$ is {\em countably generated} if there is a countably generated right $\cal B$-module ${\cal E}_0$ such that the space ${\cal E}_0 (A)$ is a dense
  subset of the space ${\cal E}(A)$ for every object $A\in \Ob ({\cal
  B})$.

The countable set $\Omega$ which generates the right $\cal A$-module
${\cal E}_0$ is referred to as a {\em generating set} for the Hilbert
$\cal A$-module $\cal E$.
\end{definition}

Note that the above definition is the same as that in \cite{Jo2,
  Mitch3} but differs from that of \cite{Mitch2}.

When the $C^\ast$-category $\cal A$ is not unital, the Hilbert $\cal
A$-modules $\Hom (-,A)_{\cal A}$ are not countably generated in
general.

\begin{definition}
A $C^\ast$-category $\cal A$ is called {\em $\sigma$-unital} if each
$C^\ast$-algebra $\Hom (A,A)_{\cal A}$ has a countable approximate
unit.
\end{definition}

When $\cal A$ is a $\sigma$-unital $C^\ast$-category, it is clear that
the Hilbert $\cal A$-modules $\Hom (-,A)_{\cal A}$ are all countably
generated.  

\begin{definition}
Let $\cal B$ be a graded $C^\ast$-category.  Then a {\em Hilbert
  $\cal B$-module} $\cal E$ is called {\em graded} if each space
  ${\cal E}(A)$ admits decompositions ${\cal E}(A) = {\cal E}(A)_0 \oplus {\cal E}(A)_1$ into vectors of degree $0$ and vectors of degree $1$ such that

\begin{itemize}

\item[$\bullet$] $\deg (\eta x) = \deg (\eta) + \deg (x)$ for all vectors $\eta \in {\cal E}(B)$ and morphisms $x\in \Hom (A,B)_{\cal B}$

\item[$\bullet$] $\deg (\langle \eta , \xi \rangle ) = \deg (\eta ) + \deg (\xi )$ for all vectors $\eta \in {\cal E}(B)$ and $\xi \in {\cal E}(A)$.

\end{itemize}

Here all addition takes place modulo $2$.
\end{definition}

\begin{definition}
Let $\cal E$ and $\cal F$ be Hilbert modules over a $C^\ast$-category, $\cal B$.  Then an {\em operator} $T\colon {\cal
  E}\rightarrow {\cal F}$ is a collection of maps $T_A \colon {\cal
  E}(A)\rightarrow {\cal F}(A)$ such that there are maps $T_A^\star
\colon {\cal F}(A)\rightarrow {\cal E}(A)$ with the property
$$\langle \eta , T_A \xi \rangle = \langle T_B^\star , \xi \rangle$$
for all vectors $\eta \in {\cal F}(B)$ and $\xi \in {\cal E}(A)$.
\end{definition}

It is shown in \cite{Mitch2} that an operator $T\colon {\cal
  E}\rightarrow {\cal F}$ is a natural
transformation, each map $T_A\colon {\cal E}(A) \rightarrow {\cal
  F}(A)$ is bounded and linear, and the collection of maps $T_A^\star$
defines an operator $T^\star$.  The operator $T^\star$ is called the
{\em adjoint} of the operator $T$.

An operator $T$ is called {\em bounded} if the norm
$$\| T \| = \sup \{ \| T_A \| \ |\ A\in \Ob ({\cal B}) \}$$
is finite.  The adjoint of a bounded operator is bounded.  

If $\cal B$ is a graded $C^\ast$-category, we write ${\cal L}({\cal
  B})$ to denote the category of all countably generated graded
  Hilbert $\cal B$-modules and bounded linear operators.  It can be
  shown (see \cite{Mitch2}) that the category ${\cal L}({\cal B})$ is
  a $C^\ast$-category.  Moreover, it is a graded $C^\ast$-category;
we define the {\em degree} of $T$ by the formula
$$\deg (T\eta ) = \deg (T) + \deg (\eta )$$

Of course, addition in the above formula takes place modulo $2$.

\begin{definition} \label{bimodule}
Let $\cal A$ and $\cal B$ be graded $C^\ast$-categories.  Then a {\em graded Hilbert $({\cal A},{\cal B})$-bimodule} is a graded $C^\ast$-functor 
${\cal E}\colon {\cal A}\rightarrow {\cal L}({\cal B})$.
\end{definition}

For each object $A\in \Ob ({\cal A})$ we write ${\cal E}(-,A)$ to
denote the corresponding graded Hilbert $\cal B$-module.  Given
another object $B\in \Ob ({\cal B})$ we have a vector space ${\cal
  E}(B,A)$.  For each morphism $x\in \Hom (A,B)_{\cal A}$ there is a bounded
operator $x\colon {\cal E}(-,A)\rightarrow {\cal E}(-,B)$.

\begin{definition}
A {\em bounded operator} $T\colon {\cal E}\rightarrow {\cal E}'$
between graded Hilbert $({\cal A},{\cal B})$-bimodules $\cal E$ and ${\cal
  E}'$ is a collection, $T$, of operators $T_A \colon {\cal
  E}(-,A)\rightarrow {\cal E}'(-,A)$ such that the norm
$$\| T \| = \sup \{ \| T_A \| \ |\ A\in \Ob ({\cal A}) \}$$
is finite.

We say that the operator $T$ has degree $k$ if each operator $T_A \colon {\cal
  E}(-,A)\rightarrow {\cal E}'(-,A)$ has degree $k$.
\end{definition}

Note that we make no assumptions here concerning naturality.  If
$T$ is a bounded operator between graded Hilbert $({\cal A},{\cal
  B})$-bimodules, and $x\in \Hom (A,A')_{\cal A}$ is a morphism in the
$C^\ast$-category $\cal A$, let us define the {\em graded commutator}
by the formula
$$[x,T] = xT_A - (-1)^{\deg (x) \deg (T)} T_{A'}x$$

The above formula only makes sense when the degree of the morphism $x$
and the operator $T$ are defined.  However, we can extend the 
definition of the graded commutator by requiring it to be linear in
each variable.

\begin{definition}
Let $\cal B$ be a $C^\ast$-category, and let $\cal E$ and $\cal F$ be
Hilbert $\cal B$-modules.  Then a {\em rank one operator} $T\colon
{\cal E}\rightarrow {\cal F}$ is an operator of the form
$$\zeta \mapsto \eta \langle \xi , \zeta \rangle$$
for elements $\eta \in {\cal F}$ and $\xi \in {\cal E}$.  We write
this operator $\eta \langle \xi , - \rangle$.  A {\em
  compact operator} is a norm-limit of finite linear combinations of rank one operators.
\end{definition}

\begin{definition}
We write ${\cal D}({\cal A},{\cal B})$ to denote the category of
graded Hilbert $({\cal A},{\cal B})$-bimodules and
bounded operators $T\colon {\cal E}\rightarrow {\cal E}'$ such that
the graded commutator $[x,T]$ is compact for all morphisms $x\in \Hom
(A,A')_{\cal A}$.  We write ${\cal K}{\cal D}({\cal A},{\cal B})$ to
denote the (non-unital) subcategory consisting of bounded operators
$T\colon {\cal E}\rightarrow {\cal E}'$ such that the composites
$xT_A$ and $T_{A'}x$ are compact for all morphisms $x\in \Hom
(A,A')_{\cal A}$.
\end{definition}

The category ${\cal D}({\cal A},{\cal B})$ is a graded $C^\ast$-category and
the subcategory ${\cal K}{\cal D}({\cal A},{\cal B})$ is a
$C^\ast$-ideal.  We can therefore form the quotient ${\cal Q}{\cal
  D}({\cal A},{\cal B}) = {\cal D}({\cal A},{\cal B}) / {\cal K}{\cal
  D}({\cal A},{\cal B})$.

\begin{definition}
Let $\cal A$ and $\cal B$ be small $\sigma$-unital graded $C^\ast$-categories.  We define the {\em $KK$-theory spectrum}
$${\mathbb K}{\mathbb K}({\cal A},{\cal B}) = \Omega {\mathbb K}{\cal
  Q}{\cal D}({\cal A},{\cal B})$$
\end{definition}

It is proved in \cite{Mitch3} that if $A$ and $B$ are
$C^\ast$-algebras, the usual $KK$-theory groups, as defined by
Kasparov in \cite{Kas2}, can be recovered as the stable homotopy
groups of the spectrum ${\mathbb K}{\mathbb K}(A,B)$.

The following results are proved in \cite{Mitch3}.

\begin{proposition}
The $C^\ast$-category ${\cal Q}{\cal D}({\cal A},{\cal B})$ is
contravariantly functorial in the variable $\cal A$ and covariantly
functorial in the variable $\cal B$.
\noproof
\end{proposition}

Hence the $KK$-theory spectrum ${\mathbb K}{\mathbb K}({\cal A},{\cal B})$ is contravariantly functorial in the variable $\cal A$ and covariantly functorial in the variable $\cal B$.

\begin{proposition}
Let $\cal B$ be a small $\sigma$-unital graded $C^\ast$-category.
Then the spectra ${\mathbb K}{\mathbb K}({\mathbb F},{\cal B})$ and
${\mathbb K}({\cal B})$ are naturally stably equivalent.  
\noproof
\end{proposition}

The main property of $KK$-theory we need in this article is a special
case of the Kasparov product.

\begin{proposition} \label{product}
Let $\cal A$ and $\cal B$ be $\sigma$-unital $C^\ast$-categories.  Then we have a
product
$${\mathbb K}({\cal A})\wedge {\mathbb K}{\mathbb K}({\cal A},{\cal
  B})\rightarrow {\mathbb K}({\cal B})$$

This product agrees with the Kasparov product when $A$ and $B$ are
$C^\ast$-algebras.  It is natural in the variable $\cal B$ in the
obvious sense, and natural in the variable $\cal A$ in the sense that
for a $C^\ast$-functor $F\colon {\cal A}\rightarrow {\cal A}'$ we
have an induced commutative diagram
$$\xymatrix{
{\mathbb K}({\cal A})\wedge {\mathbb K}{\mathbb K}({\cal A},{\cal B})
\ar[rd] \\
{\mathbb K}({\cal A})\wedge {\mathbb K}{\mathbb K}({\cal A}',{\cal B})
\ar[u]^{F^\star} \ar[d]_{F_\star} & {\mathbb K}({\cal B}) \\
{\mathbb K}({\cal A}')\wedge {\mathbb K}{\mathbb K}({\cal A}',{\cal B})
\ar[ru] \\
}$$
\noproof
\end{proposition}

Proposition \ref{fgp} gives us the following special case.

\begin{proposition} \label{moduleproduct}
Let $\cal A$ and $\cal B$ be trivially graded unital $C^\ast$-categories, and
let $\cal E$ be a finitely generated projective Hilbert $\cal A$-module.  Then we have
a natural map of spectra
$$[{\cal E}]\wedge \colon {\mathbb K}{\mathbb K}({\cal A},{\cal
  B})\rightarrow {\mathbb K}({\cal B})$$

Given a $C^\ast$-functor $F\colon {\cal A}\rightarrow {\cal A}'$ we
have a commutative diagram
$$\xymatrix{
{\mathbb K}{\mathbb K}({\cal A},{\cal B}) \ar[rd]^{[{\cal E}]\wedge} \\
& {\mathbb K}{\cal B} \\
{\mathbb K}{\mathbb K}({\cal A}',{\cal B}) \ar[uu]^{F^\star}
\ar[ru]_{F_\star [{\cal E}]\wedge } \\
}$$
\noproof
\end{proposition}

We are now ready to construct the descent map relating equivariant
$KK$-theory to crossed products.

\begin{theorem} \label{descent}
Let $A$ and $B$ be $\sigma$-unital graded $\cal G$-$C^\ast$-algebras.  Then we have a
canonical graded $C^\ast$-functor
$$D\colon {\cal Q}{\cal D}_{\cal G}(A,B)\rightarrow {\cal Q}{\cal
  D}(A\rtimes_r {\cal G}, B\rtimes_r {\cal G})$$

The $C^\ast$-functor $D$ is natural in the variables $A$ and $B$ in
the obvious sense, and natural in the variable $\cal G$ in the sense
that given a faithful functor $f\colon {\cal H}\rightarrow {\cal G}$ between
groupoids we have a commutative diagram
$$\xymatrix{
{\cal Q}{\cal D}_{\cal G}(A,B) \ar[r] \ar[d] & 
{\cal Q}{\cal D}(A\rtimes_r {\cal G}, B\rtimes_r {\cal G}) \ar[r] &
{\cal Q}{\cal D}(A\rtimes_r {\cal H}, B\rtimes_r {\cal G}) \\
{\cal Q}{\cal D}_{\cal H}(A,B) \ar[r] & 
{\cal Q}{\cal D}(A\rtimes_r {\cal H}, B\rtimes_r {\cal H}) \ar[ru] \\
}$$
\end{theorem}

\begin{proof}
Let $\cal E$ be a countably generated $\cal G$-equivariant graded Hilbert $(A,B)$-bimodule.  We
begin our construction by associating to $\cal E$ a countably
generated graded Hilbert $(A\rtimes_r
{\cal G},B\rtimes_r {\cal G})$-bimodule $D({\cal E})$.

Fix an object $a\in \Ob ({\cal G})$.  Then for each object $b\in \Ob
({\cal G})$ there is a vector space
$$D({\cal E})_0 (b,a) = \{ \eta_1 g_1 + \cdots + \eta_n g_n \ |\
\eta_i \in {\cal E}_a, g_j \in \Hom (b,a)_{\cal G} \}$$

We can define a linear contravariant functor $D({\cal E})_0 (-,a)$
from the convolution category $B{\cal G}$ to the category of vector
spaces.  The object $b\in \Ob ({\cal G})$ is mapped to the vector
space $D({\cal E})_0 (b,a)$.  The action of the category $B{\cal G}$
is defined by the formula\footnote{Note the similarity with the
  composition law in the convolution category.}
$$\left( \sum_i \eta_i g_i \right) \left( \sum_j x_j h_j \right) =
\sum_{i,j} \eta_j g_i (x_j) g_i h_j$$

There is an inner product
$$D({\cal E})_0 (c ,a)\times D({\cal E})_0(b,a) \rightarrow \Hom
(b,c)_{B\rtimes_r {\cal G}}$$
defined by the formula
$$\langle \sum_i \eta_i g_i , \sum_j \xi_j h_j \rangle = \sum_{i,j}
g_j^{-1} (\langle \eta_i , \xi_j \rangle) g_j^{-1} h_j$$

Completing the spaces $D({\cal E})_0 (b,a)$ with respect to the norms
defined by the above inner products we obtain a Hilbert $B\rtimes_r
{\cal G}$-module $D({\cal E})(-,a)$.  This Hilbert module can be
graded by saying that the sum $\sum_i \eta_i g_i$
has degree $k$ if each vector $\eta_i \in {\cal E}_a$ has degree $k$.

We can define a graded $C^\ast$-functor $D({\cal E})\colon A\rtimes_r {\cal
  G}\rightarrow {\cal L}(B\rtimes_r {\cal G})$ by mapping the object
  $a\in \Ob ({\cal G})$ to the Hilbert $B\rtimes_r {\cal G}$-module
  $D({\cal E})(-,a)$.  The action of the $C^\ast$-category $A\rtimes_r
  {\cal G}$ is defined by the formula
$$\left( \sum_i x_i g_i \right) \left( \sum_j \eta_j h_j \right) =
\sum_{i,j} x_i g_i (\eta_j ) g_i h_j$$

The $C^\ast$-functor $D({\cal E})$ is the desired graded Hilbert $(A\rtimes_r
{\cal G},B\rtimes_r {\cal G})$-bimodule associated to the graded Hilbert
$(A,B)$-bimodule $\cal E$.

Suppose we have a bounded operator $T\colon {\cal E}\rightarrow {\cal
  E}'$ between $\cal G$-equivariant graded Hilbert $(A,B)$-bimodules ${\cal
  E}$ and ${\cal E}'$.  Then we have an operator $D(T)\colon D({\cal
  E})\rightarrow D({\cal E}')$ defined by the formula
$$D(T) \left( \sum_i \eta_i g_i \right) = \sum_i (T\eta_i )g_i$$

The degree of the operator $D(T)$ is the same as that of the operator
$T$.  A straightforward calculation verifies that if $T$ is a morphism in the
category ${\cal D}_{\cal G}(A,B)$, the operator $D(T)$ is a
morphism in the category ${\cal D}(A\rtimes_r {\cal G},B\rtimes_r {\cal
  G})$.  Similarly, if the operator $T$ is a
morphism in the category ${\cal K}{\cal D}_{\cal G}(A,B)$, the operator $D(T)$ is a morphism in the
category ${\cal K}{\cal D}(A\rtimes_r {\cal G},B\rtimes_r {\cal G})$.  We thus
have a graded $C^\ast$-functor
$$D\colon {\cal Q}{\cal D}_{\cal G}(A,B)\rightarrow {\cal Q}{\cal
  D}(A\rtimes_r {\cal G}, B\rtimes_r {\cal G})$$

The desired naturality properties are easy to check.
\end{proof}

We therefore have a natural map of spectra
$$D\colon {\mathbb K}{\mathbb K}_{\cal G}(A,B)\rightarrow {\mathbb
  K}{\mathbb K}(A\rtimes_r {\cal G}, B\rtimes_r {\cal G})$$

If $G$ is a discrete group, the above map of spectra induces a map of
$KK$-theory groups $D\colon KK_G(A,B)\rightarrow
KK(A\rtimes_r G, B\rtimes_r G)$.  This
induced map is the same as the descent map defined by Kasparov in \cite{Kas3}.

A similar construction is possible if we look at full crossed
products; we do not need the details here.

\section{Assembly} \label{assembly}

Let $G$ be a discrete group, and let $X$ be a $G$-space.  Let us
assume that the group $G$ acts on the space $X$ on
the right, so the $C^\ast$-algebra $C_0 (X)$ is a
$G$-$C^\ast$-algebra, with $G$-action defined by the formula
$$(g\varphi )(x) = \varphi (xg )$$

\begin{definition}
Let $X$ be a $G$-$CW$-complex.  Then $X$ is called a {\em proper $G$-$CW$-complex} if for every point $x\in X$ the isotropy group
$$G_x = \{ x\in X \ |\ xg=x \}$$
is finite.
\end{definition}

Note that a $G$-$CW$-complex is a proper $G$-space in the usual sense
(see for example \cite{Pal}) if and only if it is a proper
$G$-$CW$-complex according to the above definition.

Recall that a $G$-space $X$ is called {\em $G$-compact} if the
quotient $X/G$ is compact.  Observe that a $G$-compact proper
$G$-$CW$-complex must be locally compact.  Any proper $G$-$CW$-complex is the direct
limit of its $G$-compact subspaces.

\begin{definition}
Let $A$ be a $\sigma$-unital $G$-$C^\ast$-algebra, and let $X$ be a proper $G$-$CW$-complex.
Then we define the {\em $G$-equivariant $K$-homology spectrum of $X$
with coefficients in $A$} to be the direct limit
$${\mathbb K}_\mathrm{hom}^G (X;A) = \lim_{\longrightarrow \atop K\
  G\mathrm{-compact}} {\mathbb K}{\mathbb K}_G (C_0 (K),A)$$
\end{definition}

Note that the $C^\ast$-algebra $C_0 (K)$ is $\sigma$-unital when $K$
is a $G$-compact proper $G$-$CW$-complex.  According to theorem \ref{descent} we have a natural map
$$D\colon {\mathbb K}{\mathbb K}_G (C_0 (K),A)\rightarrow {\mathbb
  K}{\mathbb K}(C_0 (K)\rtimes_r G , A\rtimes_r G )$$

The $C^\ast$-algebra $C_0 (K)\rtimes_r G$ is itself a finitely generated
projective Hilbert $C_0 (K)\rtimes_r G$-module.  Let us label this module ${\cal E}_K$.  Then by 
corollary \ref{moduleproduct} we have an induced map
$$[{\cal E}_K] \wedge \colon {\mathbb K}{\mathbb K}(C_0 (K)\rtimes_r G ,
A\rtimes_r G ) \rightarrow {\mathbb K}(A\rtimes_r G )$$

Composing these two maps and taking the direct limit we obtain a map
$$\beta \colon {\mathbb K}_\mathrm{hom}^G (X;A)\rightarrow {\mathbb
  K}(A\rtimes_r G )$$

\begin{definition}
The map $\beta$ is called the {\em Baum-Connes assembly map} with {\em
  coefficients} in the $G$-$C^\ast$-algebra $A$.
\end{definition}

The assembly map for the Baum-Connes conjecture with coefficients is
described in section 9 of \cite{BCH}.  The above definition is simply
a version of the standard definition at the level of spectra.  

\begin{definition}
Let $G$ be a discrete group.  Then a proper $G$-$CW$-complex $\EG$ is called
a {\em classifying space for proper actions of $G$} if for a given
subgroup $H\leq G$ the fixed point set $\EG^H$ is contractible if $H$
is finite, and empty if $H$ is infinite.
\end{definition}

Note that the classifying space $\EG$ always exists, and is unique up
to $G$-homotopy equivalence.  For details, see \cite{BCH}.

\begin{definition}
A group $G$ is said to satisfy the {\em Baum-Connes conjecture} with
{\em coefficients} in the $G$-$C^\ast$-algebra $A$ if the
Baum-Connes assembly map
$$\beta \colon {\mathbb K}_\mathrm{hom}^G (\EG ;A)\rightarrow {\mathbb
  K}(A\rtimes_r G)$$
is a stable equivalence of spectra.
\end{definition}

Again, the reader can consult \cite{BCH} for further details.  Note
that if $G$ is a finite group, we can take the classifying space $\EG$
to be a single point, and the Baum-Connes conjecture for $G$ is true
for trivial reasons.

As we explained in the introduction, the main purpose of this article
is to generalise the Baum-Connes assembly map in such a way that it
fits into the picture described by the following result.

\begin{theorem} \label{eqass}
Let $\mathbb E$ be a $G$-homotopy-invariant functor from the category
of proper $G$-$CW$-complexes to the category of spectra.  Then there
is a $G$-homotopy-invariant excisive functor ${\mathbb E}^\%$ and a
natural transformation $\alpha \colon {\mathbb E}^\% \rightarrow
{\mathbb E}$ such that the map
$$\alpha \colon {\mathbb E}^\% (G/H) \rightarrow {\mathbb E}(G/H)$$
is a stable equivalence for every finite subgroup, $H$, of the group
$G$.

Further, the pair $({\mathbb E}^\% , \alpha )$ is unique up to weak
equivalence.
\noproof
\end{theorem}

The above result is a special case of theorem 6.3 in \cite{DL}.  We
call the map $\alpha \colon {\mathbb E}^\% \rightarrow {\mathbb E}$
the {\em equivariant assembly map} associated to the functor ${\mathbb E}$.

\begin{definition}
Let $X$ be a $G$-space.  Then we write $\overline{X}$ to denote the
groupoid in which the collection of objects is that set $X$, and the
morphism sets are defined by writing
$$\Hom (x,y)_{\overline{X}} = \{ g\in G \ |\ xg = y \}$$
\end{definition}

We regard the groupoid $\overline{X}$ as a discrete groupoid; even though
there is a topology inherited from the space $X$, we do not want to
take this information into account.  If $f\colon X\rightarrow Y$ is an
equivariant map between $G$-spaces we have an induced faithful functor
$f_\star \colon \overline{X}\rightarrow \overline{Y}$.

In the literature the groupoid $\overline{X}$ is often referred to as the
{\em crossed product} of $X$ by $G$.  We do not adopt this terminology
here for fear of confusion with the other crossed products that are present.

There is an obvious natural faithful functor $i\colon \overline{X}
\rightarrow G$.  If $A$ is a $G$-$C^\ast$-algebra, then $A$ can also
be regarded as a $\overline{X}$-$C^\ast$-algebra; we associate the
$C^\ast$-algebra $A$ to each object of the groupoid $\overline{X}$,
and the morphism $g\colon A\rightarrow A$ to the element 
$g\in \Hom (x,y)_{\overline{X}} \subseteq G$.

Now let $X$ be a proper $G$-$CW$-complex.  Let $K$ be a $G$-compact subspace
of $X$.  Then we have an induced restriction map\footnote{See
  proposition \ref{restriction}.}
$$i^\star \colon {\mathbb K}{\mathbb K}_G (C_0 (K),A) \rightarrow
{\mathbb K}{\mathbb K}_{\overline{X}} (C_0 (K),A)$$

By proposition \ref{descent} there is a natural map
$$D\colon {\mathbb K}{\mathbb K}_{\overline{X}} (C_0 (K),A)\rightarrow {\mathbb
  K}{\mathbb K}(C_0 (K)\rtimes_r \overline{X} , A\rtimes_r \overline{X} )$$

\begin{definition}
Let $x\in X$.  We write ${\cal E}_K (x)$ to denote the set of collections
$$\{ \eta_y \in \Hom (x,y)_{C_0 (K)\rtimes_r \overline{X}} \ |\ y\in X \}$$
such that the formula
$$\eta_y g = \eta_z$$
is satisfied for all elements $g\in G$ such that $yg = z$.
\end{definition}

The collection of spaces ${\cal E}_K (x)$ is a Hilbert $C_0 (K)\rtimes_r
\overline{X}$-module.  The $C_0 (K)\rtimes_r \overline{X}$-action is defined by
composition of morphisms.  The inner product is defined by the formula
$$\langle \{ \eta_y \} , \{ \xi_y \} \rangle = \eta_y^\star \xi_y$$
for any point $y\in X$.

The Hilbert module ${\cal E}_K$ is not in general finitely generated and
projective.  However, we do have the following result.

\begin{proposition} \label{HilbertK}
The Hilbert $C_0 (K)\rtimes_r \overline{X}$-module ${\cal E}_K$ defines a
$K$-theory element $[{\cal E}_K] \in {\mathbb K}(C_0 (K)\rtimes_r
\overline{X})_0$.
\end{proposition}

\begin{proof}
Given a point $x\in X$, let us write $\Or (x)$ to denote the orbit of
$x$ in the space $X$ with respect to the action of the group $G$.
Observe that $\Hom (x,y)_{C_0 (K)\rtimes_r \overline{X}} = \{ 0 \}$ unless
$\Or (x) = \Or (y)$.  Let $\overline{X}|_{\Or (x)}$ denote the full
subcategory of the groupoid $\overline{X}$ in which the objects are the
points of the orbit $\Or (x)$.  The
$K$-theory spectrum ${\mathbb K}(C_0 (K)\rtimes_r \overline{X})$ is naturally
equivalent to the spectrum
$${\mathbb K}\left( \prod_{\Or (x) \in X/G} C_0 (K)\rtimes_r
  \overline{X}|_{\Or (x)} \right)$$
by proposition \ref{productcat}.

Let ${\cal E}|_{\Or (x)}$ be the Hilbert $C_0 (K)\rtimes_r \overline{X}|_{\Or
  (x)}$-module defined by restricting the Hilbert module ${\cal
  E}_K$.  Then the Hilbert module ${\cal E}|_{\Or (x)}$ is isomorphic
  to the module $\Hom (-,x)_{C_0 (K)\rtimes_r \overline{X}|_{\Or (x)}}$ and so
  is finitely generated and projective.

Hence the direct sum
$$\prod_{\Or (x) \in X/G} {\cal E}|_{\Or (x)}$$
is a finitely generated projective Hilbert module over the $C^\ast$-category
$\prod_{\Or (x) \in X/G} C_0 (K)\rtimes_r \overline{X}|_{\Or (x)}$.  It
therefore defines a $K$-theory element
$$[{\cal E}_K] \in {\mathbb K}(C_0 (K)\rtimes_r \overline{X})$$
as required.
\end{proof}

The $K$-theory element $[{\cal E}_K]$ has the following naturality property.

\begin{proposition} \label{natural1}
Let $f\colon (X,K)\rightarrow (Y,L)$ be a map of pairs of $G$-spaces,
where the subspaces $K$ and $L$ are $G$-compact.  Let
$$f^\star \colon C_0 (L)\rtimes_r \overline{X} \rightarrow C_0 (K)\rtimes_r
\overline{X} \qquad f_\star \colon C_0 (L)\rtimes_r \overline{X} \rightarrow C_0
(L)\rtimes_r \overline{Y}$$
be the obvious induced maps.  Then there is a $K$-theory element
$[\theta ] \in {\mathbb K}(C_0 (L)\rtimes_r \overline{X})$ such that $f_\star
[\theta ] = [{\cal E}_L]$ and $f^\star [\theta ]=[{\cal E}_K]$.
\end{proposition}

\begin{proof}
Let $x\in X$.  Let $\theta (x)$ to denote the set of collections
$$\{ \eta_y \in \Hom (x,y)_{C_0 (L)\rtimes_r \overline{X}} \ |\ y\in X \}$$
such that the formula
$$\eta_y g = \eta_z$$
is satisfied for all elements $g\in G$ such that $yg = z$.

The collection of spaces $\theta (x)$ is a Hilbert $C_0 (L)\rtimes_r
\overline{X}$-module.  The $C_0 (L)\rtimes_r \overline{X}$-action is defined by
composition of morphisms.  The inner product is defined by the formula
$$\langle \{ \eta_y \} , \{ \xi_y \} \rangle = \eta_y^\star \xi_y$$
for any point $y\in X$.

As in proposition \ref{HilbertK} the Hilbert $C_0 (L)\rtimes_r
\overline{X}$-module $\theta$ defines a $K$-theory element $[\theta ] \in
{\mathbb K}(C_0 (L)\rtimes_r \overline{X})$.  The formulae $f_\star
[\theta ] = [{\cal E}_L]$ and $f^\star [\theta ]=[{\cal E}_K]$ are
easy to check.
\end{proof}

Now, proposition \ref{product} gives us a map
$$[{\cal E}_K ] \wedge \colon {\mathbb K}{\mathbb K} (C_0 (K)\rtimes_r
\overline{X}, A\rtimes_r \overline{X}) \rightarrow {\mathbb K}(A\rtimes_r \overline{X})$$

We can form the composite
$$\gamma_K \colon {\mathbb K}{\mathbb K}_G (C_0 (K),A)\rightarrow
{\mathbb K}(A\rtimes_r \overline{X})$$
and take direct limits to obtain a map
$$\gamma \colon {\mathbb K}_\mathrm{hom}^G (X;A)\rightarrow {\mathbb
  K}(A\rtimes_r \overline{X})$$

\begin{lemma}
The map $\gamma \colon {\mathbb K}_\mathrm{hom}^G (X;A)\rightarrow {\mathbb
  K}(A\rtimes_r \overline{X})$ is natural for proper $G$-$CW$-complexes $X$.
\end{lemma}

\begin{proof}
Let $f\colon (X,K)\rightarrow (Y,L)$ be a map of pairs of $G$-spaces,
where the subspaces $K$ and $L$ are $G$-compact.  Then by proposition
\ref{natural1}, and naturality of the Kasparov product and descent
map, we have a commutative diagram
$$\xymatrix@=7pt{
&
{\mathbb K}{\mathbb K}_{\overline{X}} (C_0 (K),A) \ar[r] \ar[d] &
{\mathbb K}{\mathbb K}(C_0 (K)\rtimes_r \overline{X},A\rtimes_r \overline{X})
\ar[rd] \ar[d] \\
{\mathbb K}{\mathbb K}_G (C_0 (K),A) \ar[ru] \ar[d] &
{\mathbb K}{\mathbb K}_{\overline{X}} (C_0 (L),A) \ar[r] &
{\mathbb K}{\mathbb K}(C_0 (L)\rtimes_r \overline{X},A\rtimes_r \overline{X})
\ar[r] \ar[d] &
{\mathbb K}(A\rtimes_r \overline{X}) \ar[d] \\
{\mathbb K}{\mathbb K}_G (C_0 (L),A) \ar[r] \ar[ru] &
{\mathbb K}{\mathbb K}_{\overline{Y}} (C_0 (L),A) \ar[u] \ar[rd] &
{\mathbb K}{\mathbb K}(C_0 (L)\rtimes_r \overline{X},A\rtimes_r \overline{Y}) \ar[r] &
{\mathbb K}(A\rtimes_r \overline{Y}) \\
& & 
{\mathbb K}{\mathbb K}(C_0 (L)\rtimes_r \overline{Y},A\rtimes_r \overline{Y}) \ar[u]
\ar[ru] \\
}$$

Taking direct limits, the desired result follows.
\end{proof}

\begin{lemma}
Let $H$ be a finite subgroup of $G$.  Then the map
$\gamma \colon {\mathbb K}_\mathrm{hom}^G (G/H;A)\rightarrow {\mathbb
  K}(A\rtimes_r \overline{G/H})$ is a stable equivalence of spectra.
\end{lemma}

\begin{proof}
Let $i\colon H\hookrightarrow G$ be the inclusion homomorphism.  Then
by proposition \ref{gequivalence} the map $\gamma$ is equivalent to
the composite
$$\xymatrix{
{\mathbb K}{\mathbb K}_G(C_0 (G/H),A) \ar[d]^{i^\star} \\
{\mathbb K}{\mathbb K}_H(C_0 (G/H),A) \ar[d]^D \\
{\mathbb K}{\mathbb K}(C_0 (G/H)\rtimes_r H , A\rtimes_r H) \ar[d]^{[{\cal
    E}_{G/H}]\wedge} \\
{\mathbb K}(A\rtimes_r H) \\
}$$

Let $j\colon {\mathbb F}\rightarrow C_0 (G/H)$ be the inclusion
defined by writing $j (\lambda )([1]) = \lambda$ and $j (\lambda
)(x)=0$ if $x\neq [1]$.  Let $+$ denote the one point topological
space.  Then ${\mathbb F} = C_0 (+)$ and we have a commutative diagram:
$$\xymatrix{
{\mathbb K}{\mathbb K}_G(C_0 (G/H),A) \ar[d]^{i^\star} \\
{\mathbb K}{\mathbb K}_H(C_0 (G/H),A) \ar[d]^D \ar[r]^{j^\star} &
{\mathbb K}{\mathbb K}_H(C_0 (+), A) \ar[d]^D \\
{\mathbb K}{\mathbb K}(C_0 (G/H)\rtimes_r H , A\rtimes_r H) \ar[d]^{[{\cal
    E}_{G/H}]\wedge} \ar[r]^{j^\star} & {\mathbb K}{\mathbb K}(C_0 (+)\rtimes_r H ,
A\rtimes_r H) \ar[d]^{[{\cal E}_{+}]\wedge} \\
{\mathbb K}(A\rtimes_r H) \ar[r]^1 & {\mathbb K}(A\rtimes_r H) \\
}$$

A straightforward calculation tells us that the composite $j^\star
i^\star \colon {\mathbb K}{\mathbb K}_G (C_0 (G/H),A)\rightarrow
{\mathbb K}{\mathbb K}_H (C_0 (+),A)$ is a stable equivalence of
spectra.  The composite map on the right, $\beta \colon {\mathbb
  K}{\mathbb K}_H (C_0 (+),A)\rightarrow {\mathbb K}(A\rtimes_r H)$,
is the Baum-Connes assembly map.  Since the group $H$ is finite, the
space $+$ is a model for the
classifying space $\underline{E}H$ and the map $\beta$ is a stable
equivalence of spectra.  Hence the map
$$\gamma \colon {\mathbb K}_\mathrm{hom}^G (G/H;A)\rightarrow {\mathbb
  K}(A\rtimes_r \overline{G/H})$$
is a stable equivalence as claimed.
\end{proof}

The canonical functor $i\colon \overline{X}\rightarrow G$
gives us an induced map
$$i_\star \colon {\mathbb K}(A\rtimes_r \overline{X}) \rightarrow {\mathbb
  K}(A\rtimes_r G)$$
by proposition \ref{crossfunctor}

\begin{lemma}
The composite map $i_\star \gamma \colon {\mathbb
  K}_\mathrm{hom}^G(X;A)\rightarrow {\mathbb K}(A\rtimes_r G)$ is the
  Baum-Connes assembly map.
\end{lemma}

\begin{proof}
Let $K$ be a $G$-compact subspace of $X$.  The naturality properties
of the various descent maps and products give us a commutative diagram
$$\xymatrix{
{\mathbb K}{\mathbb K}_G (C_0 (K),A) \ar[r]^-D \ar[dd]_{i^\star} &
{\mathbb K}{\mathbb K}(C_0 (K)\rtimes_r G,A\rtimes_r G) \ar[r]^-{[{\cal
    E}_K]\wedge} \ar[d]_{i^\star} & {\mathbb K}(A\rtimes_r G) \\
& {\mathbb K}{\mathbb K} (C_0 (K)\rtimes_r \overline{X} ,A\rtimes_r G)
\ar[ru]_-{[{\cal E}_K]\wedge} \\
{\mathbb K}{\mathbb K}_{\overline{X}} (C_0 (K),A) \ar[r]^-D &
{\mathbb K}{\mathbb K}(C_0 (K)\rtimes_r \overline{X},A\rtimes_r \overline{X})
\ar[u]^{i_\star} \ar[r]^-{[{\cal E}_K]\wedge} & {\mathbb K}(A\rtimes_r \overline{X}) \ar[uu]_{i_\star} \\
}$$

Here the top row is the Baum-Connes assembly map, and the
composite $[{\cal E}_K]\wedge \circ D\circ i^\star$ is the map
$\gamma$.  Taking direct limits, the desired result follows.
\end{proof}

Observe that the equivariant $K$-homology functor ${\mathbb
  K}^G_\mathrm{hom}(-,A)$ is $G$-homotopy-invariant and excisive.  We
  can therefore use the above three lemmas to apply theorem \ref{eqass} to
the study of the Baum-Connes assembly map.  We immediately obtain the
following result.

\begin{theorem} \label{main}
Let ${\mathbb E}^\%$ be a $G$-homotopy-invariant excisive functor from
the category of proper $G$-$CW$-complexes to the category of spectra.
Suppose we have a natural transformation $\alpha \colon {\mathbb E}^\%
(X)\rightarrow {\mathbb K}(A\rtimes_r \overline{X})$ such that the map
$$\alpha \colon {\mathbb E}^\% (G/H) \rightarrow {\mathbb K}(A\rtimes_r \overline{G/H})$$
is a stable equivalence for every finite subgroup, $H$, of the group
$G$.

Then up to stable equivalence the composite $\alpha i_\star \colon
{\mathbb E}^\% (X)\rightarrow {\mathbb K}(A\rtimes_r G)$ is the Baum-Connes assembly map.
\noproof
\end{theorem}

We conclude by using the above result to give an alternative
description of the Baum-Connes assembly map.  To formulate
it, we need to introduce one more piece of machinery from \cite{DL}.

\begin{definition}
Let $G$ be a discrete group.  Then we define the {\em orbit category},
$\Or (G)$, to be the category in which the objects are $G$-spaces,
$G/H$, where $H$ is a subgroup of $G$, and the morphisms are
$G$-equivariant maps.
\end{definition}

An {\em $\Or (G)$-spectrum} is a functor from the category $\Or (G)$
to the category of symmetric spectra.
Our main example of an $\Or (G)$-spectrum is defined by writing
$${\mathbb E}(G/H) = {\mathbb K}(A\rtimes_r \overline{G/H})$$
where $A$ is a given $G$-$C^\ast$-algebra.

The following result is proved in \cite{DL}.

\begin{theorem} \label{DL2}
Let $\mathbb E$ be an $\Or (G)$-spectrum.  Then there is a
$G$-homotopy-invariant excisive functor, ${\mathbb E}^\%$, from the
category of $G$-$CW$-complexes to the category of spectra such that ${\mathbb
  E}^\% (G/H) = {\mathbb E}(G/H)$ whenever $H$ is a subgroup of $G$.

Further, given a functor $\mathbb F$ from the category of
$G$-$CW$-complexes to the category of spectra, there is a
natural transformation
$$\beta \colon ({\mathbb F}|_{\Or (G)})^\% \rightarrow {\mathbb F}$$
such that the map
$$\beta \colon ({\mathbb F}|_{\Or (G)})^\% (G/H) \rightarrow {\mathbb F}(G/H)$$
is a stable equivalence for every subgroup, $H$, of the group $G$.
\noproof
\end{theorem}

\begin{theorem}
Consider the $\Or (G)$-spectrum
$${\mathbb E}(G/H) = {\mathbb K}(A\rtimes_r \overline{G/H})$$

Let $X$ be a path-connected space, and let $c\colon X \rightarrow +$ be the constant map.  Then up to
stable equivalence the induced map
$$c_\star \colon {\mathbb E}^\% (X) \rightarrow {\mathbb E}^\% (+)$$
is the Baum-Connes assembly map.
\end{theorem}

\begin{proof}
Consider the functor
$${\mathbb F} \colon X\mapsto {\mathbb K}(A\rtimes_r \overline{X})$$

Then ${\mathbb F}|_{\Or (G)} = {\mathbb E}$.  Let $c\colon X\rightarrow \pi_0
(X)$ be the map defined by sending a point of a topological space $X$ to its path-component.  
By theorem \ref{DL2} we have a commutative diagram
$$\xymatrix{
{\mathbb E}^\% (X) \ar[r]^\beta \ar[d]_{c_\star} & {\mathbb E}(X) \ar[d]^{c_\star} \\
{\mathbb E}^\% \pi_0 (X) \ar[r]^\beta & {\mathbb E}\pi_0 (X) \\
}$$
where the map $\beta \colon {\mathbb E}^\% (G/H)\rightarrow {\mathbb
  E}(G/H)$ is a stable equivalence whenever $H$ is a subgroup of $G$.

If the space $X$ is path-connected, by theorem \ref{main} the composite
$c_\star \beta$ is equivalent to the Baum-Connes assembly map.
Observe that the map $\beta \colon {\mathbb E}^\% (+)\rightarrow
{\mathbb E}(+)$ is a stable equivalence.\footnote{Since the quotient
  $G/G$ is just the one-point topological space.}

Hence the map $c_\star \colon {\mathbb E}^\% (X)\rightarrow
{\mathbb E}^\%(+)$ is equivalent to the Baum-Connes assembly map as claimed.
\end{proof}

In particular, it follows from the above result that the map
$$c_\star \colon {\mathbb E}^\% (\EG ) \rightarrow {\mathbb E}^\% (+)$$
is the Baum-Connes assembly map.  This fact is already stated in
\cite{DL} for the Baum-Connes assembly map {\em without} coefficients.
An alternative proof of the above result, at least in the
coefficient-free case, can be found in \cite{HP}.

In the coefficient-free case, the above theorem is used in
\cite{Lueck2} to explicitly calculate the rationalisations of the
$K$-theory groups $K_n C^\star_r (G)$ when $G$ is a discrete group
that satisfies the Baum-Connes conjecture.  Taking coefficients into
account, a similar calculation should be possible to calculate the
rationalisations of the groups $K_n (A\rtimes_r G)$.

\bibliographystyle{plain}

\bibliography{data.bib}

\end{document}